\documentclass[12pt,reqno]{amsart}
\usepackage{latexsym}

\newcommand{\K}{{\mathbf K}}
\renewcommand{\epsilon}{\varepsilon}

\renewcommand{\Re}{{\operatorname{Re}}}
\renewcommand{\Im}{{\operatorname{Im}}}
\newcommand{\tr}{{\operatorname{Tr}}}
\newcommand{\sym}{{\operatorname{Sym}}}

\newcommand{\HS}{{\operatorname{HS}}}
\newcommand{\crit}{{\operatorname {crit}}}

\newcommand{\szego}{Szeg\"o }
\newcommand{\nhat}{\raisebox{2pt}{$\wh{\ }$}}

\newcommand{\inv}{^{-1}}
\newcommand{\kahler}{K\"ahler }
\newcommand{\sqrtn}{\sqrt{N}}
\newcommand{\wt}{\widetilde}
\newcommand{\wh}{\widehat}
\newcommand{\PP}{{\mathbb P}}
\newcommand{\R}{{\mathbb R}}
\newcommand{\C}{{\mathbb C}}

\newcommand{\CP}{\C\PP}
\renewcommand{\d}{\partial}
\newcommand{\dbar}{\bar\partial}
\newcommand{\ddbar}{\partial\dbar}
\newcommand{\U}{{\rm U}}

\newcommand{\half}{{\frac{1}{2}}}
\newcommand{\vol}{{\operatorname{Vol}}}

\newcommand{\SU}{{\operatorname{SU}}}
\newcommand{\FS}{{{\operatorname{FS}}}}

\renewcommand{\phi}{\varphi}

\newcommand{\ccal}{\mathcal{C}}

\newcommand{\hcal}{\mathcal{H}}
\newcommand{\ical}{\mathcal{I}}
\newcommand{\kcal}{\mathcal{K}}

\newcommand{\lcal}{\mathcal{L}}
\newcommand{\ncal}{\mathcal{N}}

\newcommand{\rcal}{\mathcal{R}}
\newcommand{\ocal}{\mathcal{O}}
\newcommand{\scal}{\mathcal{S}}

\newcommand{\jcal}{\mathcal{J}}
\newcommand{\al}{\alpha}
\newcommand{\be}{\beta}
\newcommand{\ga}{\gamma}

\newcommand{\La}{\Lambda}
\newcommand{\la}{\lambda}
\newcommand{\ep}{\varepsilon}
\newcommand{\de}{\delta}
\newcommand{\De}{\Delta}
\newcommand{\om}{\omega}
\newcommand{\Om}{\Omega}

\newtheorem{theo}{{\sc Theorem}}[section]
\newtheorem{cor}[theo]{{\sc Corollary}}
\newtheorem{conj}[theo]{{\sc Conjecture}}

\newtheorem{lem}[theo]{{\sc Lemma}}
\newtheorem{prop}[theo]{{\sc Proposition}}

\newenvironment{rem}{\medskip\noindent{\it Remark:\/} }{\medskip}
\newenvironment{defin}{\medskip\noindent{\it Definition:\/} }{\medskip}
\newenvironment{claim}{\medskip\noindent{\it Claim:\/} }{\medskip}

\title[Critical points and supersymmetric vacua II]
{Critical points and supersymmetric vacua, II: Asymptotics and
extremal metrics}

\author{Michael R. Douglas}
\address{NHETC and Department of Physics and Astronomy,
Rutgers University,
Piscataway, NJ 08855--0849, USA;
and I.H.E.S., Bures-sur-Yvette, France}
\email{mrd@physics.rutgers.edu}

\author{Bernard Shiffman}
\address{Department of Mathematics, Johns Hopkins University, Baltimore,
MD 21218, USA} \email{shiffman@math.jhu.edu}

\author{Steve Zelditch}
\address{Department of Mathematics, Johns Hopkins University, Baltimore,
MD 21218, USA} \email{zelditch@math.jhu.edu}

\thanks{Research partially supported by
DOE grant DE-FG02-96ER40959  (first author) and NSF grants
DMS-0100474 (second author) and DMS-0302518 (third author).}

\begin{document}

\begin{abstract}
 Motivated by the vacuum selection problem of string/M theory,
we study a  new geometric invariant of a
 positive Hermitian line bundle $(L, h)
\to M $ over a compact \kahler manifold: the expected distribution
of critical points  of a Gaussian random holomorphic section $s
\in H^0(M, L)$ with respect to the Chern connection $\nabla_h$. It
is a measure on $M$ whose total mass is the average number
$\ncal^\crit_h$ of critical points of a random holomorphic
section. We are interested in the metric dependence of
$\ncal^\crit_h$, especially metrics $h$ which minimize
$\ncal^\crit_h$. We concentrate on the asymptotic minimization
problem for the sequence of  tensor powers $(L^N, h^N) \to M$ of
the line bundle and their critical point densities
$\kcal^\crit_{N,h}(z)$. We prove that $\kcal^\crit_{N,h}(z)$ has a
complete asymptotic expansion in $N$ whose coefficients are
curvature invariants of $h$. The first two terms in the expansion
of $\ncal^\crit_{N,h}$ are topological invariants of $(L, M)$. The
third term is a topological invariant plus a constant $\beta_2(m)$
(depending only on the dimension $m$ of $M$) times the Calabi
functional $\int_M \rho^2 d\vol_h$, where $\rho$ is the scalar
curvature of the \kahler metric $\om_h:=\frac i2\Theta_h$. We give an
integral formula for $\beta_2(m)$ and show, by a computer assisted
calculation,  that $\beta_2(m) > 0$ for $m \leq 5$, hence that
$\ncal^\crit_{N,h}$ is asymptotically minimized  by the Calabi
extremal metric (when one exists).  We conjecture that $\beta_2(m)
> 0$ in all dimensions, i.e.\ the Calabi extremal metric is
always the asymptotic minimizer.

\end{abstract}

\maketitle

\tableofcontents

\section{Introduction}  This paper is the second in a series of
articles \cite{DSZ, DSZ2}  on the statistics of vacua in string/M
theory and associated  effective supergravity theories.
Mathematically, vacua are critical points $\nabla s(z) = 0$
 of a holomorphic section $s \in H^0(M, L)$ of
a line bundle $L \to M$ over a complex manifold relative to a
connection $\nabla,$ which we always choose to be the Chern
connection $\nabla_h$ of a Hermitian metric $h$ on $L$.
Equivalently, they are critical points of the norm $|s(z)|_h$ of
$s$ relative to $h$. Our motivation to study critical points of
holomorphic sections in this metric sense comes initially from
physics, where critical points of this kind model  extremal black
holes in addition to string/M  vacua (cf.\ \cite{FGK, St, M}). But
we also find the statistics of critical points to have an
independent geometric interest. A basic statistical quantity, the
average number of critical points of a random holomorphic section,
defines a new geometric invariant of a positive Hermitian
holomorphic line bundle, and in this paper we show that its
asymptotic minima are given by Calabi extremal metrics.

 The physical setting for  vacua of string/M theory (and for extremal black holes)
 is a holomorphic line bundle   over the moduli space of complex
 structures on a  polarized  Calabi-Yau manifold.
Supersymmetric vacua are critical points
 of a holomorphic section (known as a superpotential) relative to
 the Weil-Petersson metric. The program of studying
 statistics of critical points of Gaussian   random
holomorphic sections was proposed by the first author  \cite{AD,
D, DD} as a means of dealing with the large number of string/M
vacua. There exists at this time no reasonable selection principle
to decide which superpotential nor which of its critical points
gives the vacuum state which correctly describes our universe in
string/M theory, so it   makes sense to study the statistics of
 vacua of random superpotentials.

In our first paper \cite{DSZ}, we gave explicit formulas for the
expected distribution  of critical points of sections of general
holomorphic line bundles over any complex manifold, including
those which arise in physics. The formulas, recalled in \S
\ref{background} (cf.\ Theorem \ref{KNcrit2}), involve complicated
complex symmetric matrix integrals, and it is difficult to see how
the expected distribution  of critical points depends on the
metric $h$. The purpose of this article is therefore to study a
purely geometric simplification of the physical problem where $(L,
h)$ is a positive Hermitian line bundle over a (usually compact)
\kahler manifold $M$ and where the Gaussian measure on $H^0(M, L)$
is derived from the inner product induced by $h$. Our aim is to
understand the metric dependence of the statistics  of the random
critical point set
\begin{equation}  \label{Critdef} Crit(s, h) = \{z: \nabla_h (s) = 0 \}
\end{equation}
   of a
Gaussian random section of $H^0(M, L)$. We note that $Crit(s, h)
\cup Z_s =  \{z: d \;  |s(z)|_h^2 = 0\}$ where
$Z_s$ is the zero set of $s$.

From the probabilistic viewpoint, the critical points of random
holomorphic sections relative to the Chern connection $\nabla_h$
of a fixed Hermitian metric on $L$ determine a {\it point process}
on $M$, that is, a measure on the configuration space of
finite subsets of $M$, which gives the probability density of a
finite subset being  the critical point set of a holomorphic
section. The process  is determined by its $n$-point correlation
functions $\K^\crit_{n} (z_1, \dots, z_n)$, which give the
probability density of critical points occurring simultaneously
at the points
$z_1, \dots, z_n \in M$.

In this article, we focus on the $1$-point correlation,
namely the expected distribution of critical points
\begin{equation*}  \K^\crit_{h}  =\int
 \bigg[\,\sum_{z\in Crit (s, h)} \delta_{z}\bigg]\,d\ga_h(s)
\end{equation*} where $\delta_{z}$ is the Dirac point mass at $z$,
and where ${\gamma_h}$ is the Gaussian measure probability $\gamma_h$
on
$H^0(M, L)$ induced by
$h$ and the \kahler form $\om_h:=\frac i2 \Theta_h$ (see
\S\ref{background}). We showed in
\cite{DSZ} that
$\K^\crit_{h}$ is a smooth form on $M$ if $L$ is sufficiently
positive. In particular, we are interested in the expected
(average) number of critical points
\begin{equation}\ncal^\crit_{h}=\int_M
\K^\crit_{h}= \int\#Crit(s,h)\,d\ga_h(s) \end{equation} of a random
section, a purely geometric invariant of a positive Hermitian
holomorphic line bundle $(L, h) \to M$. The methods of this paper also
give results on general $n$-point correlations between critical points
and their scaling asymptotics in the sense of \cite{BSZ1}, but we do
not carry out the analysis here.

Since it is a crucial point, let us explain why   the number $\#
Crit(s,h)$ is a (non-constant) random variable on $H^0(M,
L)$, unlike the number of zeros of $m$ independent sections which
is a topological invariant of $L$.  As indicated above, connection
critical points  are the same as critical points of $|s(z)|_h^2$
for which $s(z) \not= 0$, or equivalently as critical points of
$\log |s(z)|_h$ (see \cite{DSZ} for the simple proof). Hence,
there are critical points of each Morse index $\geq m$ (see
\cite{B, DSZ}), and only the alternating sum of the number of
critical points of each index is a topological invariant.  Another
way to understand the metric dependence of the number of critical
points is to write the covariant derivative in a local frame $e_L$
as
\begin{equation}
\nabla_{z_j} s=\left( \frac{\d f}{\d z_j} -f \frac{\d K}{\d
z_j}\right)e_L= e^{K}\frac{\d }{\d z_j}\left(e^{-K}\,
f\right)e_L\;,\qquad \nabla_{\bar z_j} s= \frac{\d f}{\d\bar
z_j}\, e_L\;,
\end{equation}
where we locally express a section as $s=fe_L$, and $K=-\log
\|e_L\|^2_h$. Hence, the critical point equation
\begin{equation}\label{CPE}  \left( \frac{\d f}{\d z_j} -f \frac{\d K}{\d
z_j}\right) = 0 \end{equation} in the local frame fails to be
holomorphic when the connection form is only smooth.

Although $Crit(s, h)$ and $\# Crit(s, h)$ depend on $h$, it is not
 clear at the outset whether $\ncal^\crit_h$ is a topological
invariant or to what degree it depends on the metric $h$. To
investigate the metric dependence of $\K^\crit_h$ and
$\ncal^\crit_h$ we consider their asymptotic behavior as we take
powers $L^N$ of $L$. As in \cite{SZ, BSZ1}, it is natural to
expect that the density and number of critical points will have
 asymptotic expansions which reveal their metric dependence.

We therefore let  $\K^\crit_{N,h}(z)$  denote the expected
distribution of critical points of random holomorphic sections
$s_N \in H^0(M, L^N)$ with respect to the Chern connection and
Hermitian Gaussian measure induced by $h$, as given by
\eqref{inner}--\eqref{gauss} in \S \ref{background}.  We also let
\begin{equation}\ncal^\crit_{N,h}=\int_M \K^\crit_{N,h}(z)
\end{equation} denote the expected number of critical points.
The covariant  derivative associated to $h^N$ has the
semi-classical form
\begin{equation}
\nabla_{z_j} s_N=\left( \frac{\d f}{\d z_j} -Nf \frac{\d K}{\d
z_j}\right)e_L^{\otimes N}= e^{NK}\frac{\d }{\d
z_j}\left(e^{-NK}\, f\right)e_L^{\otimes N}\;,\qquad \nabla_{\bar
z_j} s_N= \frac{\d f}{\d\bar z_j}\, e_L^{\otimes N}\;.
\end{equation}

We let $\Theta_h=\ddbar K$
denote the curvature form of $(L,h)$.
Our first result is a complete asymptotic expansion for the
expected distribution of critical points for powers $L^N \to M$ in terms
of curvature invariants of the \kahler metric $\om_h=\frac
i2\Theta_h$\;:  

\begin{theo} \label{UNIVDIST}  For any positive Hermitian  
line bundle $(L, h) \to (M,
\omega_h)$ over any compact \kahler manifold with $\om_h = \frac
i2\Theta_h$, the expected critical point distribution of
holomorphic sections of $L^N$ relative to the Hermitian Gaussian
measure induced by $h$ and $\om_h$ has an asymptotic expansion of the
form
$$N^{-m}\,\K_{N,h}^\crit(z)\sim \{b_0  + b_1(z) N\inv  +b_2(z)
N^{-2} + \cdots\} \frac{\omega_h^m}{m!}\;,$$
 where the coefficients
$b_j=b_j(m)$ are curvature invariants of order $j$ of $\omega_h$. In
particular, $b_0$ is the universal constant
\begin{equation}\label{b0}b_0 = \pi^{-{m+2\choose 2}}  \int_{\sym(m,\C)
\times \C} \left|\det(2H H {}^*-|x|^2I)\right|\,e^{- \langle (H,
x), (H,x) \rangle} \,dH\,dx\,,\end{equation}
  $b_1=\be_1 \rho$, where $\rho$ is the  scalar curvature of
$\om_h$ and $\be_1$ is a universal constant, and $b_2 $ is a
quadratic curvature polynomial.
\end{theo}

Here, $\sym(m,\C)$ is the space of $m \times m$ complex symmetric
matrices.  It follows that the density of critical points is
asymptotically uniform relative to the curvature volume form  with
a universal asymptotic density $b_0(m)$. The values of the
constant $b_0$ for low dimensions are:
\begin{eqnarray*} &\displaystyle b_0(1)=\frac5{3\,\pi}\;,\quad
 b_0(2)=\frac{59\cdot 2}{3^3\,\pi^2}
\;,\quad b_0(3)=\frac{637\cdot 3!}{3^5\,\pi^3}\;,\quad
b_0(4)=\frac{6571\cdot 4!}{3^7\,\pi^4}\;,\\ &\displaystyle
b_0(5)=\frac{65941\cdot 5!}{3^9\,\pi^5}\;,\quad
b_0(6)=\frac{649659\cdot 6!}{3^{11}\,\pi^6}\;.\end{eqnarray*} (See
Appendix 1.)

With only minor changes in the proofs, our methods give
refinements of the asymptotic results which take the Morse indices
of the critical points into account. By the Morse index $q$ of a
critical point, we mean its Morse index as a critical point of
$\log \|s\|_{h^N}$; it is well known that $m\le q\le 2m$ for
positive line bundles \cite{B}. Thus we let
$\K^\crit_{N,q,h}(z)=\K^\crit_{N,q}(z)$ denote the expected
distribution of critical points of $\log \|s_N\|_{h^N}$ of Morse
index $q$,   and we let $\ncal^\crit_{N,q,h}$ denote the expected
number of these critical points. Thus we have
\begin{equation} \K^\crit_{N,h}(z)=\sum_{q=m}^{2m}\K^\crit_{N,q,h}(z)\,,\qquad
\ncal^\crit_{N,h}=\sum_{q=m}^{2m}\ncal^\crit_{N,q,h}\;.\end{equation}
We obtain a similar asymptotic expansion for the distribution of
critical points of given Morse index:
\begin{theo} \label{univmorse}  Let $M,L,h,\om_h,\K^\crit_{N,q,h}$ be
as above. Then the expected density of critical points of
Morse index $q$ of random sections in $H^0(M,L^N)$ is given by
$$N^{-m}\,\K_{N,q,h}^\crit(z)\sim \{b_{0q}  + b_{1q}(z) N\inv  + b_{2q}(z)
N^{-2} + \cdots\}\frac {\om_h^m}{m!}\;,\qquad m\le q\le 2m\;,$$ where
the $b_{jq}=b_{jq}(m)$ are curvature invariants of order $j$ of
$\omega_h$. In particular, $b_{0q}$ is given by the integral in
\eqref{b0} with the domain of integration $\sym(m, \C)\times \C$
replaced by
\begin{equation}\label{Smk}{\bf S}_{m,k}:=\{(H,x)\in
\sym(m,\C)\times\C: \mbox{\rm index}(2HH^*-|x|^2I)
=k\}\;,\end{equation} with $k=q-m$.
\end{theo}

 The coefficient integrals are very complicated to
evaluate except in dimension one, where we obtain the following
formula:

\begin{theo} \label{critRS} Let $(L,h)$ be a positive line bundle on a compact complex
curve $C$ of genus $g$.  Then the expected number of saddle points and of
local maxima, respectively, of random holomorphic sections of $L^N$ are
given by:
\begin{eqnarray*}\ncal^\crit_{N,1,h}&=& \frac 43 \,c_1(L)\, N + \frac 89\,
( 2g-2) +
\left(\frac 1{27\pi}\,
\int_C \rho^2 \om_h\right)N\inv +O(N^{-2})\;,\\ \ncal^\crit_{N,2,h}&=&
\frac 13\,c_1(L)\, N - \frac 19\, ( 2g-2) +
\left(\frac 1{27\pi}\,
\int_C \rho^2 \om_h\right)N\inv +O(N^{-2})\;,\end{eqnarray*} where
$\om_h=
\frac i2 \Theta_h$, and $\rho$ is the Gaussian curvature of the
metric $\om_h$.\end{theo}

Hence, the expected total number of critical points in dimension 1 is
$$\ncal^\crit_{N,h}= \frac 53
\,c_1(L)\, N + \frac 79\, ( 2g-2) +
\left(\frac 2{27\pi}\,
\int_C \rho^2 \om_h\right)N\inv +O(N^{-2})\;.$$ In dimension 1, the
coefficients in Theorem \ref{univmorse} are:
$$b_{01}(1)=\frac 4{3\pi}\,, \ b_{02}(1)=\frac 1{3\pi}\,,\
b_{11}(1)= -\frac 8{9\pi}\,\rho, \  b_{12}(1)= \frac
1{9\pi}\,\rho, \  b_{21}(1)= b_{22}(1)= \frac
1{27\,\pi}\,\rho^2.$$ In Appendix 1, we provide a table of
numerical values of $b_{0q}(m)$ (for $m\le 6$), which we computed
using Maple.

In \S \ref{s-exact}, we study the case where  $M=\CP{^m}$
and $(L,h)$ is the hyperplane section bundle $\ocal(1)$
with the Fubini-Study metric. We show that in this case the expected number
$\ncal^\crit_{N,q,h}$ of critical points of  Morse
index $q$ of sections in $H^0(\CP^m,\ocal(N))$ is a rational function of
$N$, and we provide an exact formula (Proposition
\ref{numcrit}) for
this expected number. In Appendix~1, we use Proposition
\ref{numcrit} and  a computer assisted
computation to give explicit formulas for $\ncal^\crit_{N,q,h}$ in
dimensions
$\le 4$.  For example, in dimension 1, the expected total number of critical
points, as computed in
\cite{DSZ}, is
$$\ncal^\crit_{N,h}(\CP^1) =\frac{5N^2-8N+4}{3N-2},$$ and in dimension 2,
the  expected total number is
$$\ncal^\crit_{N,h}(\CP^2)=
{\frac {59\,{N}^{5}-231\,{N}^{4}+375\,{N}^{3}-310\,{N}^{2}+132\,N-24}{
 \left( 3\,N-2 \right) ^{3}}}\;.$$

Thus we gain a quantitative sense of how many more critical points
there are in the metric sense in comparison with the
classical critical point equation $\d f  = 0$. In the case of $\ocal(N)
\to
\CP^1$, whose sections are polynomials of degree $N$, we may view this
classical critical point equation as a connection critical point
equation by viewing  the derivative $\frac{\partial }{\partial z}$ as a
flat meromorphic connection with pole at $\infty$. Alternately, it is
the Chern connection of a singular Hermitian metric. The critical
point equation being purely holomorphic, the number of critical
points of a generic section is a constant $N-1$.  All critical
points relative to this connection are saddle points. By
comparison, sections $s_N \in H^0(\CP^1, \ocal(N))$
have, on average, an additional
$\sim
\frac{N}{3}$ local maxima and $\sim \frac{N}{3}$ additional
saddles relative to the Fubini-Study connection. The study of critical
points relative to meromorphic connections (known as Minkowski vacua)
is simpler than that relative to Chern connections and will be
explored further in a subsequent work.

As  a corollary of Theorem \ref{univmorse}, we find the rather
surprising fact that
 the asymptotics of the
expected number of critical points is a topological invariant of
the bundle $L \to M$  to two orders in $N$.

\begin{cor} \label{EN}  Let $(L,h)\to (M,\om_h)$ be a
positive holomorphic line bundle on a compact \kahler manifold,
with $\om_h= \frac i2 \Theta_h$. Then the expected  number of
critical points of Morse index $q$ ($m\le q\le 2m$) of random
sections in $H^0(M,L^N)$ has the asymptotic expansion
\begin{eqnarray}\ncal^\crit_{N,q,h} &\sim
&\left[\frac{\pi^m b_{0q}}{m!}\,c_1(L)^m\right]N^m +
\left[\frac{\pi^m\be_{1q}}{(m-1)!}\, c_1(M)\cdot
c_1(L)^{m-1}\right]N^{m-1} \nonumber
\\&+& \left[\be_{2q}\int_M\rho^2d\vol_h + \be'_{2q}\, c_1(M)^2\cdot
c_1(L)^{m-2} + \be''_{2q}\, c_2(M)\cdot c_1(L)^{m-2}\right]
N^{m-2}+\cdots\,,\nonumber\end{eqnarray} where
$b_{0q},\be_{1q},\be_{2q},\be'_{2q},\be''_{2q}$ are universal constants
depending only on the dimension $m$.
\end{cor}

We recall that the alternating sum of the $\ncal^\crit_{N,q,h}$ is
topological:
$$\sum  (-1)^{m+q}\,\ncal^\crit_{N,q,h} = c_n(L^N \otimes
T^{*1,0}) = c_1(L)^m \,N^m + \cdots\,,$$ and hence by Corollary \ref{EN},
\begin{equation}\label{alt}\sum_{q=m}^{2m} (-1)^{m+q}\,
b_{0q}(m)=\frac{m!}{\pi^m}\;.\end{equation} Since each
$b_{0 q}$ is strictly positive and their sum equals
$b_0$, it follows from \eqref{alt} that
\begin{equation}\label{lower} b_{0}(m)> \frac{m!}{\pi^m}\;.\end{equation}

Corollary \ref{EN} shows that the metric dependence of
$\ncal^\crit(h^N)$ is a `lower-order' effect, and hence renews the
question  to what degree $ \ncal^\crit_{q,h}$, or at least  the
asymptotic expansion of $ \ncal^\crit_{N,q,h}$,  is a topological
invariant. We see from Corollary \ref{EN} that the expansion is
not topological provided that the constant $\be_{2q}=\be_{2q}(m)$
does not vanish.  Computations in dimensions $\le 5$ show that
$\be_{2q}$ is positive in these dimensions (see Corollary
\ref{LESS4}), and we  conjecture that $\be_{2q}(m)>0$ for all $m$.
As we now explain, this conjecture is suggested by a connection
between extremals of $\ncal^\crit_{N, h}$ and Calabi extremal
metrics
  \cite{Ca1, Ca2, T, Don}.

   This connection involves a notion of asymptotic minimality of
   $\ncal_{N,h}^\crit$. To introduce it, we revisit the motivating  problem of determining how $\ncal^\crit_h$ varies
as $h$ varies over Hermitian metrics on $L$. One could consider
all Hermitian metrics on $L$, but we focus on the smaller class of
positively curved Hermitian  metrics,
$$P(M, L) = \{h: \frac{i}{2} \Theta_h \;\; \mbox{is a positive}\;\;
(1,1)\mbox{-form} \;\}.$$   If we fix one such metric $h_0 =
e^{-K_0}$, the others may be expressed as $h_{\phi}:= e^{\phi}
h_0$ with $\phi \in C^{\infty}(M).$  It is reasonable to
conjecture that $\ncal^\crit_{h_{\phi}}$ is unbounded as
$h_{\phi}$ varies over $P(M, L)$,  since by  (\ref{CPE}) the
number of critical points of a section should be `large' if the
`degree' of the connection form $-\partial K$ is `large'. Here, $K
= -\log h = K_0 - \phi$ in a local frame. On the other hand,
$\ncal^\crit_h$ is bounded below by  $|c_m(L \otimes T^{* 1,0})|$,
and it is plausible that it has a
 minimum. It would be interesting to determine this minimal
metric (assuming one exists and is unique), but it is difficult to
solve the critical point equation $\delta \ncal^\crit_h = 0$.

We therefore consider the simpler asymptotic problem. Since the
first two leading coefficients in the expansion of
$\ncal_{N,h}^\crit$ are topological, we try to find  metrics for
which the first non-topological term is critical. The first
non-topological term is  a multiple of the Calabi functional
$$\int_M \rho_h^2 \, dV_h\;,$$ where $\rho_h$ is the scalar curvature
 of the \kahler metric $\om_h= \frac i2 \Theta_h$, and
$dV_h = \frac 1{m!}\om_h^m$. Thus the problem of finding
metrics which are critical for the metric invariant
$\ncal_{N,h}^\crit $ is  closely related to the problem of finding
critical points of Calabi's functional. In keeping with our
intuition that $\ncal_{N,h}^\crit$ should have a minimizer, we
note that critical points of Calabi's functional are necessarily
minima (cf.\ \cite{Ca1, Ca2, Hw}).

Existence of critical metrics  is one of the fundamental problems
in complex geometry, and we refer to \cite{Don, T,  Y2} for
background. It was suggested by S.-T. Yau \cite{Y1, Y3} that
existence of a canonical metric should be  related to the
stability of $M$. One class of canonical metrics are Hermitian
metrics $h$ for which $\Theta_h$ is a \kahler metric of constant
scalar curvature, i.e. for which $\rho_h$ is constant. By a
theorem due to S. Donaldson \cite[Cor.~5]{Don}, there exists at
most one \kahler metric of constant scalar curvature in the
cohomology class of $ \pi\, c_1(L)$. Hence if there exists such a
metric of constant scalar curvature, there exists a unique
Hermitian metric minimizing Calabi's functional.

This leads us to make the following definition:

\begin{defin} Let $L \to M$ be an ample holomorphic line bundle
over a compact \kahler manifold.   For $h \in P(M,L)$ and $ m\le q\le 2m$,
we say that
$\ncal_{N,q,h}^\crit$ (resp.\ $\ncal_{N,h}^\crit$) is {\it
asymptotically minimal\/} if for all $ h_1
\neq h$ in $P(M,L)$, there exists $ N_0=N_0(h_1)$ such that
\begin{equation}\ncal^\crit_{N,q,h} < \ncal^\crit_{N,q,h_1}\quad
(\mbox{resp. } \ncal^\crit_{N,h} < \ncal^\crit_{N,h_1})\qquad
\mbox{for }\  N \geq N_0 \; .  \end{equation}
\end{defin}

Assuming $(M, L)$ has a Hermitian metric $h$ minimizing Calabi's
functional, we see from Corollary \ref{EN} that  $\ncal_{N,q,h}^\crit$
(resp.\ $\ncal_{N,h}^\crit$) is   asymptotically minimal as long as
$\be_{2q}(m)> 0$ (resp.\ $\be_2(m)>0$). Since we believe this to be the
case for all dimensions, we state the following conjecture.

\begin{conj} \label{conj}  Let $h
\in P(M,L)$, $m\le q\le 2m$.  Then the following are
equivalent:
\begin{itemize} \item $\ncal_{N,h}^\crit$ is asymptotically minimal,
\item $\ncal_{N,q,h}^\crit$ is asymptotically minimal (and hence
$\ncal_{N,q',h}^\crit$  is asymptotically minimal $\forall q'$),
\item $h$
minimizes Calabi's functional.
\end{itemize}
\end{conj}

In Lemma \ref{K0}, we show that
\begin{equation} \label{K0INTRO} \be_{2q}(m)= \frac
1{4\,\pi^{m+2\choose 2}} \int_{{\bf S}'_{m,q-m}} \ga(H) \,
\left|\det(2HH^*-|x|^2I) \right| \,e^{ -\langle (H, x),(H, x)
\rangle}\, dH\, dx\;,
\end{equation}  where
$$\gamma(H) = \frac{1}{2} |H_{11}|^4 - 2 |H_{11}|^2 +
1.$$  It is unfortunately  difficult to determine  from
this formula whether $\be_{2q}(m)$ is non-zero and (if so) what sign
it has.  In \S \ref{s-coeff}, we transform
\eqref{K0INTRO} to a rather complicated, but more elementary,
integral  (Lemma \ref{LI4}), which can be evaluated by a routine
(but long) calculation. For low dimensions, we perform this
calculation using Maple  to obtain:

\begin{theo} \label{ASYMPMIN}  The constants $\be_{2q}(m)$
are positive  for $m\le 5$, and hence
Conjecture \ref{conj} is true for  $\dim M \leq 5$.
\end{theo}

 The calculation in dimension $5$ was done by B.
Baugher \cite{BB}.  In particular, we have:

\begin{cor} \label{LESS4} Suppose that  $\dim M \leq 5$
and that $L$ possesses a metric $h$ for which the scalar curvature of
$\om_h=\frac i2\Theta_h$ is constant.  Then $h$ is the unique
metric on $L$ such that $\ncal_{N,h}^\crit$ is asymptotically
minimal.\label{calabi}\end{cor}

Thus, for instance, the Fubini-Study metric $h$ on the hyperplane
section bundle $\ocal(1) \to \CP^m$ is the unique metric on $\ocal(1)$
such that $\ncal_{N,h}^\crit$ is asymptotically minimal, at least for
$m\le 5$.

 Additional evidence that $\beta_{2q}(m) > 0$ was found by
B. Baugher \cite{BB}. By studying patterns in computer-assisted
calculations of the terms arising from the computation of the
integral in  (\ref{K0INTRO}),  Baugher discovered an identity for
$\beta_{2q}(m)$  for dimensions $m \leq 5$ which  he
conjectures  is valid in all dimensions. Baugher
 showed that this identity implies that $\beta_{2q}(m)
> 0$ in all dimensions, and hence Baugher's conjecture
implies  Conjecture \ref{conj}. Baugher's conjecture and identity
are stated in Appendix 2.

We close the introduction with some comments on the organization
of the paper. The proof of Theorems \ref{UNIVDIST} and
\ref{univmorse} is based on the Tian-Yau-Zelditch asymptotic
expansion of the \szego kernel $\Pi_N(z, w)$ \cite{Ze, T1, Y2, Lu}
and on formulas from \cite{DSZ} for the density of critical
points. We then need to evaluate the coefficients explicitly to
obtain concrete results linking geometry to numbers of critical
points. Once we know the leading coefficient is universal, we may
calculate it for $\ocal(N) \to \CP^m$ and this is done in \S
\ref{s-exact} and Appendix~1. Unfortunately, the Fubini-Study
metric is not useful for finding the sign of $\be_2$ since it is
impossible to separate out the topologically invariant terms from
the Calabi functional for this metric. Hence in \S \ref{s-coeff},
we analyze  instead the case of $M=\CP^1 \times E^{m-1}$ for $E$
an elliptic curve, where the relevant topological terms vanish.
This leads to an explicit integral which we analyze by a variant
of the Itzykson-Zuber formula in random matrix integrals.

We thank X. X. Chen and Z. Lu for discussions of Calabi's
functional and references to the literature.

\section{Background}\label{background}

Let  $(L, h)  \to M$ be a Hermitian holomorphic line bundle over a
complex  manifold $M$, and let $\nabla = \nabla_h$ be its Chern
connection, i.e.\ the unique connection of type $(1,0)$
on $L$  compatible with both the metric and complex
structure of $L$. Thus, it satisfies $\nabla'' s = 0$ for any
holomorphic section $s$ where $\nabla = \nabla' + \nabla''$ is the
splitting of the connection into its $L \otimes T^{*1 ,0}$, resp.\
$L \otimes T^{*0, 1}$ parts. It follows that
\begin{equation} Crit(s, h) = \{z: \nabla_h' s (z) = 0 \}.
\end{equation}
 We denote by $\Theta_h=d\d\log h = -\ddbar \log
h$ the curvature of $h$ and $\omega_h = \frac{i}{2} \Theta_h.$

We now introduce the Gaussian measures $\gamma_h$, called {\it
Hermitian Gaussian measures} in \cite{DSZ} which we use
exclusively in this paper. They are determined by the inner
product
\begin{equation}\label{inner}\langle s_1, s_2 \rangle = \int_M h(s_1(z),
s_2(z)) dV_h(z)\end{equation} on $H^0(M, L)$, where $dV_h= \frac
1{m!}\om_h^m$.  By definition,
\begin{equation}\label{gauss}d\gamma_{h}(s)=\frac 1
{\pi^d}e^{-\|c\|^2} dc\;,\qquad  s=\sum _{j=1}^d
c_je_j,\end{equation} where $dc$ is Lebesgue measure and $\{e_j\}$
is an orthonormal basis for $\scal$ relative to $\langle,
\rangle$.

\begin{defin} \label{DEFK} The {\it expected distribution of critical
points\/} of $s \in
\scal \subset H^0(M, L)$ with respect to  $\gamma_h$  is defined
by
\begin{equation}  \K^\crit_{h}  =\int_{H^0(M, L)}
 \bigg[\,\sum_{z\in Crit (s, h)} \delta_{z}\bigg]\,d\ga_h(s)
\end{equation} where $\delta_{z}$ is the Dirac point mass at $z$;
i.e.,
\begin{equation}  \big( \K^\crit_{h}, \phi\big) =
\int_{H^0(M, L)} \left[\sum_{z: \nabla_h s(z) = 0} \phi(z)
\right]\, d\gamma_{h}(s)\,, \qquad \phi\in \ccal^\infty(M)\,.
\end{equation} The density of
$\K^\crit_h$ with respect to $dV_h$ is denoted
$\kcal_h^\crit(z)$; i.e.,
$$\K^\crit_h=\kcal_h^\crit(z)\,dV_h\;.$$
\end{defin}

\subsection{Formulas for the expected distribution of critical
points}\label{s-formulas}
 Let $(L, h) \to (M,\om)$ be a Hermitian
holomorphic line bundle on an $m$-dimensional \kahler manifold. We
say that  $ H^0(M, L)$ has the $2$-jet spanning property if  all
possible values and derivatives of order $\le 2$ are attained by
the global sections $s \in H^0(M, L)$ at every point of $M$. In
\cite{DSZ}  we showed that if $ H^0(M, L)$ has the $2$-jet spanning
property, then $\K^\crit_h$ is absolutely continuous with respect to
$dV_h$, and we obtained an integral formula for $\kcal^\crit_{h} (z_0)$ in
terms of the
\szego kernel $\Pi(z,\bar w)$  of $H^0(M, L) $ with respect to $h$.  To
describe this formula, we choose normal coordinates about $z_0\in M$ and
define the following matrices:
\begin{eqnarray}
A(z_0)&=& \big(\nabla_{z_j}
\nabla_{\bar w_{j'}}\Pi\big)\,,\label{A}\\
B(z_0)&=&  \Big[\big(\tau_{j'q'}\,\nabla_{z_j}\nabla_{\bar
w_{q'}}\nabla_{\bar w_{j'} }\Pi\big) \quad
\big(\nabla_{z_j}\Pi\big)
\Big] \,,\label{B}\\[8pt]
C(z_0)&=&\left[
\begin{array}{cc}
\big(\tau_{jq}\tau_{j'q'}\, \nabla_{z_q}\nabla_{z_j}\nabla_{\bar
w_{q'}}\nabla_{\bar w_{j'}}\Pi\big)&
\big(\tau_{jq}\, \nabla_{z_q}\nabla_{z_j}\Pi\big) \\[8pt]
\big(\tau_{j'q'}\,\nabla_{\bar w_{q'}}\nabla_{\bar
w_{j'}}\Pi\big)&
\Pi\end{array}\right] \,,\label{C}\\
&&\tau_{jq}=\sqrt 2 \ \ \mbox{ if }\ j<q\;,\quad
\tau_{jj}=1\;,\label{tau}\\
&& \qquad \qquad 1\le j\le m\,,\ 1\le j\le q\le m\,,\ 1\le j'\le
q'\le m\,,\nonumber
\end{eqnarray}
where the \szego kernel $\Pi=\Pi(z,w)$ and
its derivatives are evaluated at $(z_0,0;z_0,0)$ (see \S
\ref{s-szego} and \eqref{szegolift}). In the above matrices, $j,q$
index the rows, and
$j',q'$ index the columns.  Note that $A,B,C$ are $m\times m,\,m\times
d_m,\, d_m\times d_m$ matrices, respectively, where
$$d_m= \dim_\C(\sym(m,\C)\times\C) = \frac{m^2+m+2}2\;.$$ We then let
\begin{equation}
\La(z_0) = C(z_0)-B(z_0)^*A(z_0)^{-1}B(z_0) \;.\label{Lambda}
\end{equation}

The matrices $A,B,C$ give the second moments of the  joint
probability distribution of the random variables $\nabla s(z_0)$
and $\nabla^2 s (z_0)$ on $\scal$.

\begin{theo} \cite{DSZ} \label{KNcrit2} Let $(L, h) \to M$ denote a positive
holomorphic line bundle with the 2-jet spanning property.  Give
$M$ the volume form $dV_h= \frac 1{m!} \left(\frac
i2\Theta_h\right)^m$ induced from the curvature of $L$. Then the
expected density of critical points relative to $dV_h$ is given by
$$\kcal^\crit_{h} (z) =
\frac{\pi^{-{m+2\choose 2}}}{\det A(z) \det\La(z)}
\int_{\sym(m,\C) \times \C} \left|\det(HH{}^*-|x|^2I)\right|e^{
-{\langle \La(z)^{-1}(H, x),(H, x) \rangle}}\,dH\,dx\,. $$
\end{theo}
Here, $H \in \sym(m, \C)$ is a complex symmetric matrix, $dH$ and
$dx$ denote Lebesgue measure, and $\Lambda\inv$ is the Hermitian
operator on the complex vector space $\sym(m, \C) \times \C$
described as follows:

Let $S^{jq},\ 1\le j \le q \le m$, be the basis for $\sym(m, \C)$
given by  $$(S^{jq})_{j'q'} = \left\{\begin{array}{ll}\frac
1{\sqrt 2}(\de_{jj'}\de_{qq'} + \de_{qj'}\de_{jq'}) & \mbox{ for }\
j<q\\ \de_{jj'}\de_{qq'} & \mbox{ for }\ j=q.\end{array}\right.$$
I.e., for $j<q$, $S^{jq}$ is the matrix with $\frac 1{\sqrt 2}$ in
the $jq$ and $qj$ places and 0 elsewhere, while $S^{jj}$ is the
matrix with 1 in the $jj$ place and 0 elsewhere. We note that
$\{S^{jq}\}$ is an orthonormal basis (over $\C$) for $\sym(m,\C)$
with respect to the Hilbert-Schmidt Hermitian inner product
\begin{equation} \label{HS} \langle S,T\rangle_\HS = \tr
(ST^*).\end{equation}  For $H =(H_{jq})\in\sym(m,\C)$, we have
\begin{equation}\label{Hhat}H=\sum_{1\le j\le q\le m} \wh H_{jq}\,E^{jq}\;, \qquad \wh
H_{jq}=\tau_{jq} H_{jq}\;,\end{equation} where $\tau_{jq}$ is
given by \eqref{tau}. Lebesgue measure $dH$ (with respect to the
Hilbert-Schmidt norm) is given by
$$dH=\prod_{j\le q}d\Re \wh H_{jq}\wedge d\Im \wh H_{jq}$$

Writing
$$\Lambda=\left[
\begin{array}{cc}
\big(\Lambda_{jq}^{j'q'}\big)& \big(\Lambda_{jq}^{0}\big)
 \\[8pt]\big(\Lambda_{0}^{j'q'}\big) &\Lambda_0^0
\end{array}\right]\;,$$
we then define
\begin{equation}\label{precisedef}
\langle \La(z)^{-1}(H, x),(H, x) \rangle = \sum
(\La\inv)_{jq}^{j'q'}\wh H_{jq}\overline{\wh H_{j'q'}}
+2\Re\sum(\La\inv)_{jq}^0\wh H_{jq} \bar x +
(\La\inv)_0^0|x|^2.\end{equation}

To study the asymptotics, we consider powers $L^N$ and we use the
following result.
\begin{cor}  \label{KNcrit}
With the same notation and assumptions as above,  the density of
the  expected distribution $\K^\crit_{N,h}$ of critical points
of random sections $s_N\in H^0(M,L^N)$ relative to $dV_h$ is given
by
$$\kcal^\crit_{N,h} (z) =
\frac{\pi^{-{m+2\choose 2}}}{\det A_N(z) \det\La_N(z)}
\int_{\sym(m,\C) \times \C} \left|\det(HH{}^*-|x|^2I)\right|e^{
-{\langle \La_N(z)^{-1}(H, x),(H, x) \rangle}}\,dH\,dx\,.$$ where
\begin{equation}
\La_N(z_0) = C_N(z_0)-B_N(z_0)^*A_N(z_0)^{-1}B_N(z_0)
\;,\label{LambdaN}
\end{equation}
\begin{eqnarray}
A_N(z_0)&=& \Big[\big(\nabla_{z_j}\nabla_{\bar
w_{j'}}\Pi_N\big)\Big]\,,\label{AN}\\ B_N(z_0)&=&
\Big[\big(\tau_{j'q'}\,\nabla_{z_j}\nabla_{\bar
w_{q'}}\nabla_{\bar w_{j'} }\Pi_N\big) \quad
\big(N\,\nabla_{z_j}\Pi_N\big)
\Big] \,,\label{BN}\\[8pt]
C_N(z_0)&=&\left[
\begin{array}{cc}
\big(\tau_{jq}\tau_{j'q'}\, \nabla_{z_q}\nabla_{z_j}\nabla_{\bar
w_{q'}}\nabla_{\bar w_{j'}}\Pi_N\big)&
\big(\tau_{jq}N\, \nabla_{z_q}\nabla_{z_j}\Pi_N\big) \\[8pt]
\big(\tau_{j'q'}N\,\nabla_{\bar w_{q'}}\nabla_{\bar
w_{j'}}\Pi_N\big)&
N^2\,\Pi_N\end{array}\right] \,,\label{CN}\\
&&\tau_{jq}=\sqrt 2 \ \ \mbox{ if }\ j<q\;,\quad
\tau_{jj}=1\;,\nonumber\\
&& \qquad \qquad 1\le j\le m\,,\ 1\le j\le q\le m\,,\ 1\le j'\le
q'\le m\,,\nonumber
\end{eqnarray}
where $\Pi_N$ and its derivatives are evaluated at
$(z_0,0;z_0,0)$.
\end{cor}

\begin{proof} Rescale  $z_j = \wt z_j/\sqrtn$.  Then the curvature
of $L^N$ is given by
$$\Theta_{h^N} =N\Theta_h=\frac 12 \sum d\wt z_j\wedge d\overline{\wt
z_j},$$ so that the $\wt z_j$ are normal coordinates (at a point
$z_0$) for the curvature of $L^N$.  Apply Theorem~\ref{KNcrit2},
using the coordinates $\{\wt z_j\}$ to obtain $ A, B,  C,  \La$ .
Since $d\wt V = N^m\,dV$ and the transformation $(A,\La)\mapsto
(NA,N^2\La)$  introduces a factor $N^{-m}$, we let $A_N=N A,\
B_N=N^{3/2} B,\ C_N= N^2 C$ to obtain the desired
formula.\end{proof}

We also have a formula for the density of critical points of specific
Morse indices:

\begin{theo} \label{Morse} Under the above assumptions,  the density relative
to $dV_h$ of the  expected distribution $\K^\crit_{N,q, h}$ of critical
points  of Morse index $q$  of
$\log\|s_N \|_h$ for  random sections $s_N\in H^0(M,L^N)$  is given by
$$\kcal^\crit_{N,q,h} (z) =
\frac{\pi^{-{m+2\choose 2}}}{\det A_N(z) \det\La_N(z)} \int_{{\bf
S}_{m,q-m}}
\left|\det(HH{}^*-|x|^2I)\right|e^{ -{\langle
\La_N(z)^{-1}(H, x),(H, x) \rangle}}\,dH\,dx\,. $$ where
$${\bf S}_{m,k}=\{(H,x)\in \sym(m,\C)\times\C: \mbox{\rm index}(HH^*-|x|^2I)
=k\}\;.$$
\end{theo}

\begin{proof} The case $N=1$ is given as Theorem~6 in \cite{DSZ}. The general
case follows immediately by rescaling as in the proof of Corollary
\ref{KNcrit}.\end{proof}

Recall that the index of a nonsingular Hermitian matrix is the
number of its negative eigenvalues, and the Morse index of a
nondegenerate critical point of a real-valued function is the
index of its (real) Hessian.

\subsection{The \szego kernel}\label{s-szego}

As in our previous work, it is useful to lift the analysis
on positive line bundles $L \to M$ to the associated principal
$S^1$ bundle $X \to M$.  Sections then become scalar functions and
it is simpler to formulate various asymptotic properties for
powers $L^N$ \cite{BSZ1, BSZ2}. The same analysis is also
useful for general line bundles although the asymptotic results
no longer hold.

Given a holomorphic line bundle $L$ and a Hermitian metric $h$ on
$L$, we obtain a Hermitian metric $h^*$ on the dual line bundle
$L^*$ and we define the associated circle bundle by  $X=\{\la \in
L^* : \|\la\|_{h^*}= 1\}$.   Thus, $X$ is the boundary of the disc
bundle $D = \{\la \in L^* : \rho(\la)>0\}$, where
$\rho(\la)=1-\|\la\|^2_{h^*}$. When $(L, h)$ is a positive line
bundle, the  disc bundle $D$ is strictly pseudoconvex in $L^*$,
hence $X$ inherits the structure of a strictly pseudoconvex CR
manifold.  When $L$ is negative, as is the case for the line
bundles relevant to string theory, $X$ is pseudo-concave. We endow
$X$ with the contact  form $\al=
-i\partial\rho|_X=i\dbar\rho|_X$ and the associated
 volume form
\begin{equation}\label{dvx}dV_X=\frac{1}{m!}\al\wedge
(d\al)^m=\al\wedge\pi^*dV_M\,.\end{equation}

We  define the Hardy space $\hcal^2(X) \subset \lcal^2(X)$ of
square-integrable CR functions on $X$, i.e., functions that are
annihilated by the Cauchy-Riemann operator $\dbar_b$  and are
$\lcal^2$ with respect to the inner product
\begin{equation}\label{unitary} \langle  F_1, F_2\rangle
=\frac{1}{2\pi}\int_X F_1\overline{F_2}dV_X\,,\quad
F_1,F_2\in\lcal^2(X)\,.\end{equation}   We let $r_{\theta}x
=e^{i\theta} x$ ($x\in X$) denote the $S^1$ action on $X$ and
denote its infinitesimal generator by $\frac{\partial}{\partial
\theta}$. The $S^1$ action on $X$ commutes with
$\bar{\partial}_b$; hence $\hcal^2(X) = \bigoplus_{N =0}^{\infty}
\hcal^2_N(X)$ where $\hcal^2_N(X) = \{ F \in \hcal^2(X):
F(r_{\theta}x) = e^{i N \theta} F(x) \}$. A section $s_N$ of $L^N$
determines an equivariant function $\hat{s}_N$ on $L^*$ by the
rule
$$\hat{s}_N(\lambda) = \left( \lambda^{\otimes N}, s_N(z)
\right)\,,\quad \la\in L^*_z\,,\ z\in M\,,$$ where
$\lambda^{\otimes N} = \lambda \otimes \cdots\otimes \lambda$. We
henceforth restrict $\hat{s}$ to $X$ and then the equivariance
property takes the form $\hat s_N(r_\theta x) = e^{iN\theta} \hat
s_N(x)$. The map $s\mapsto \hat{s}$ is a unitary equivalence
between $H^0(M, L^{ N})$ and $\hcal^2_N(X)$.

We let $e_L$ be a nonvanishing local section, or local frame, of $L$.
 As above, we write
\begin{equation}\label{Kpot}\|e_L(z)\|_h^2 = e^{-K(z,
\bar z)}\;.\end{equation}
Thus, a positive line bundle $L$ induces  the \kahler
form
$\om_h=\frac i2 \ddbar K$ with \kahler potential $K$.

The \szego kernel  $\Pi_N(x,y)$ is the kernel of the  orthogonal projection $\Pi_N :
\lcal^2(X)\rightarrow
\hcal^2_N(X)$; it is defined by
\begin{equation} \Pi_N F(x) = \int_X \Pi_N(x,y) F(y) dV_X (y)\,,
\quad F\in\lcal^2(X)\,. \label{PiNF}\end{equation}

Let $\{s^N_j=f_j e_L^{\otimes N}:j=1,\dots,d_N\}$ be an
orthonormal basis for $H^0(M,L^N)$. Then $\{\hat s_j^N\}$ is an orthonormal basis
of $\hcal^2(X)$, and the \szego kernel can be written in the form
\begin{equation}\label{szego} \Pi_N(x,y)=\sum_{j=1}^{d_N}\hat s^N_j(x)
\overline{\hat s^N_j(y)}\;.\end{equation}  It is the lift of the
section \begin{equation}\label{FNdef}
\wt \Pi_N(z,\bar w): =
F_N (z,\bar w)\,e_L^{\otimes N}(z) \otimes\overline {e_L^{\otimes
N}(w)}\,,\end{equation} where
\begin{equation}\label{FN}F_N(z,\bar w)=
\sum_{j=1}^{d_N}f_j(z) \overline{f_j(w)}\;.\end{equation}

We let $(z,\theta)$ denote the coordinates of the point
$x=e^{i\theta}\|e_L(z)\|_he_L^*(z)\in X$. The equivariant lift of a section
$s=fe_L^{\otimes N}\in H^0(M,L^N)$ is given explicitly by
\begin{equation}\label{lift}\hat s(z,\theta) =
e^{iN\theta} \|e_L^{\otimes N}\|_h
f(z) = e^{N\left[-\half K(z,\bar z) +i\theta \right]} f(z)\;.\end{equation}
The \szego kernel is then given by
\begin{equation}\label{szegolift}\Pi_N(z,\theta;w,\phi) =
e^{N\left[-\half K(z,\bar z)-\half K(w,\bar w)
+i(\theta-\phi)\right]}F_N(z,\bar w)\;.\end{equation}

\subsubsection{The connection} We denote by $H = \ker \alpha$ and obtain a splitting $T_X
= H \oplus \C \frac{\partial}{\partial
\theta}$ into horizontal and vertical spaces. The Chern connection
$\nabla$ on $L^N$ then lifts to  $X$ as the horizontal derivative $d^H$, i.e.
\begin{equation}\label{horiz}
(\nabla s_N)\nhat = d^H\hat s_N\,.\end{equation}

To describe the connection explicitly, we choose
local holomorphic coordinates $\{z_1,\dots,z_m\}$ in $M$, and we write
$$\nabla=\nabla' +\nabla'',\qquad \nabla' s_N =
\sum dz_j\otimes\nabla_{z_j}s_N,\qquad\nabla'' s_N =
\sum d\bar z_j\otimes\nabla_{\bar z_j}s_N\;.$$
In particular,  $(\nabla'' s_N)\nhat = \dbar_b \hat s_N$, which vanishes when the
section
$s_N$ is holomorphic, or equivalently, when $\hat s_N\in \hcal_N^2(X)$.

For a section $s_N=fe_L^{\otimes N}$ of
$L^N$, we have
\begin{equation}\label{covar}
\nabla_{z_j} s_N=\left( \frac{\d f}{\d z_j} -Nf \frac{\d K}{\d
z_j}\right)e_L^{\otimes N}= e^{NK}\frac{\d }{\d
z_j}\left(e^{-NK}\, f\right)e_L^{\otimes N}\;,\qquad \nabla_{\bar
z_j} s_N= \frac{\d f}{\d\bar z_j}\, e_L^{\otimes N}\;.
\end{equation}

We also write
\begin{equation} d^H\Pi_N(z,\theta;w,\phi)=
\sum dz_j\otimes\nabla_{z_j}\Pi_N+
\sum d\bar w_j\otimes\nabla_{\bar w_j}\Pi_N\;,\end{equation}
where $d^H$ is the horizontal derivative on $X\times X$. (We used the fact that the
horizontal derivatives of $\Pi_N$ with respect to the $\bar z_j$ and $w_j$ variables
vanish.)   By
\eqref{szegolift}--\eqref{covar}, we have
\begin{eqnarray}\nabla_{z_j}\Pi_N &=& e^{N\left[-\half K(z,\bar z)-\half K(w,\bar w)
+i(\theta-\phi)\right]}
\left( \frac{\d }{\d z_j} -N \frac{\d K }{\d z_j}(z,\bar
z)\right)F_N(z,\bar w)\;,\label{nablaz}\\
\nabla_{\bar w_j}\Pi_N &=& e^{N\left[-\half K(z,\bar z)-\half K(w,\bar w)
+i(\theta-\phi)\right]}
\left( \frac{\d }{\d \bar w_j} -N \frac{\d K }{\d \bar w_j}(w,\bar
w)\right)F_N(z,\bar w)\;.\label{nablaw}
\end{eqnarray}

\section{Alternate formulas for the density of critical
points}\label{s-alternate}

The integrals in Theorems
\ref{KNcrit2}--\ref{Morse} are difficult to evaluate because of the
absolute value sign, which prevents a direct application of Wick
methods.  To compute the densities, we shall replace our integral by
another one which can be evaluated by
residue calculus in certain cases. This new integral is given by the
following lemma:
\begin{lem}\label{IZint} Let $\La$ be a positive definite Hermitian
operator  on
$\sym(m,\C)\times \C$. Then
\begin{eqnarray*}&& \hspace{-.2in}\frac{1}{\pi^{d_m}\,\det\La}
\int_{\sym(m,\C) \times \C} \left|\det(HH{}^*-|x|^2I)\right|e^{
-{\langle \La^{-1}(H, x),(H, x) \rangle}}\,dH\,dx\\&&
=\ \frac{(-i)^{m(m-1)/2}}{(2\pi)^m\prod_{j=1}^mj!} \,\lim_{\ep' \to
0^+}\lim_{\ep \to
0^+}
\int_{\R^m } \int_{\R^m} \int_{\U(m)} \!\! \frac{\Delta(\xi)\,
\Delta(\lambda)\; |\prod_j \la_j| \,e^{i\langle \xi, \lambda
\rangle}  e^{- \epsilon |\xi|^2 -\epsilon' |\lambda|^2}}{\det\left[
i\wh D(\xi)\rho(g)\La\rho(g)^* + I\right] }\, dg\, d \xi\, d\lambda,
\end{eqnarray*}
 where
\begin{itemize}
\item
$\Delta(\lambda) = \Pi_{i < j} (\lambda_i - \lambda_j)$, \item
$dg$ is unit mass Haar measure on $\U(m)$, \item $\wh D(\xi)$ is
the Hermitian operator on $\sym(m,\C)\oplus \C$ given by
$$ \wh D(\xi) \big((H_{jk}),x\big) =
\left( \left(\frac{\xi_j+\xi_k}{2}\,H_{jk}\right),\
-\left(\textstyle\sum_{q=1}^m \xi_q\right) x \right)\;,$$ \item
 $\rho$ is the representation  of $\U(m)$ on
$\sym(m,\C)\oplus \C$ given by
$$\rho(g)(H, x) = (gHg^t,x)\;.$$
\end{itemize}
\end{lem}
 The integrand is analytic in $\xi, g$ but rather complicated.
 Its principal features are:
\medskip

\begin{itemize}

\item $\Delta(\xi), \Delta(\lambda)$ are homogeneous polynomials
of degree $m (m - 1)/2$, and  $ |\prod_j \la_j|$ is homogeneous of
degree $m$;

\item $P_{g, z} (\xi) = \det\left[ i\wh
D(\xi)\rho(g)\La(z)\rho(g)^* + I\right] $ is a (family of)
polynomial(s) in $\xi$ of degree $m(m + 1)/2 + 1$ with no real
zeros $\xi \in \R^m$.  But the polynomials are not elliptic (or
even hypo-elliptic), that is, do not satisfy $|P(\xi)| \geq C
|\xi|^{m(m + 1)/2 + 1}$ (or any other power $|\xi|^{\mu}$).
Indeed, for large $|\xi|$ we may drop the second term $I$ and find
that the growth at infinity is that of $\det\left[ i\wh D(\xi)
\right]$. Since $\wh{D}(\xi)$ is a diagonal matrix as described in
Theorem \ref{I}, its determinant is a product of linear
polynomials in $\xi$, and hence vanishes along a union of real
hyperplanes.

\item The ratio $p_g(\xi) =  \frac{\Delta(\xi)}{\det\left[
i\wh D(\xi)\rho(g)\La(z)\rho(g)^* + I\right] }$
is thus  a rational function in $\xi$ which is  a `symbol' of
order $-m - 1$, i.e. each $\xi$-derivative decays to one extra
order. Repeated partial integrations in $d\xi$   using $\frac{1}{1
+ |\lambda|^2} [I -
 \Delta_{\xi}]
e^{i \langle \lambda, \xi \rangle} = e^{i \langle \lambda, \xi
\rangle}$ simultaneously lowers the order in both $\xi$ and
$\lambda$ by two and renders the $d \lambda$ integral absolutely
convergent without the Gaussian factor.
\end{itemize}

The proof of Lemma \ref{IZint} is given in \S\ref{SIMPLIFY} below.

As a consequence, we
have the following alternative formula for the expected critical point
density:
\begin{theo} \label{I} Under the hypotheses of
Theorem~\ref{KNcrit2} and notation of Lemma \ref{IZint}, the density of the
expected distribution of critical points of sections of $H^0(M,L^N)$ is
also given by:
\begin{equation*}
\kcal^\crit_{N,h}(z) =  \frac{c_m}{\det
A_N}\lim_{\ep' \to 0^+}  \int_{\R^m } \lim_{\ep \to 0^+}\int_{\R^m}
\int_{\U(m)}
\frac{\Delta(\xi)\,
\Delta(\lambda)\; |\prod_j \la_j| \,e^{i\langle \xi, \lambda
\rangle}   e^{- \epsilon |\xi|^2 -\epsilon'
|\lambda|^2}}{\det\left[ i\wh D(\xi)\rho(g)\La_N(z)\rho(g)^* +
I\right] }\, dg\, d \xi\, d\lambda,
\end{equation*}
 where
$$c_m=\frac{(-i)^{m(m-1)/2}}{2^m\,\pi^{2m}\,
\prod_{j=1}^mj!}\;.$$\end{theo}

\begin{proof}  The formula follows by  combining  Corollary
\ref{KNcrit} and Lemma \ref{IZint}.  \end{proof}

In \S\ref{s-exact} we shall use
 Theorem \ref{I} to calculate the density of
critical points for random sections $s_N\in H^0(\CP^m,\ocal(N))$
of the $N$-th power of the hyperplane bundle.  In this case the
$\U(m)$ integral drops out, and the integral can be evaluated as an
iterated integral without the Gaussian factor $e^{- \epsilon
|\xi|^2 -\ep' |\lambda|^2}$.

We  also have an alternative formula for the Morse index
densities, which follows by a similar argument (given in \S
\ref{Iq-proof}):

\begin{theo} \label{Iq} Under the above assumptions, the  density of the
expected distribution of critical points of Morse index q of $\log\|s_N \|_h$
is also given by:
\begin{equation*}
\kcal^\crit_{N,q,h}(z) =  \frac{m!\,c_m}{\det
A_N}\lim_{\ep' \to 0^+}  \int_{Y_{2m-q} } \lim_{\ep \to 0^+}\int_{\R^m}
\int_{\U(m)}
\frac{\Delta(\xi)\,
\Delta(\lambda)\; |\prod_j \la_j| \,e^{i\langle \xi, \lambda
\rangle}   e^{- \epsilon |\xi|^2 -\epsilon'
|\lambda|^2}}{\det\left[ i\wh D(\xi)\rho(g)\La_N(z)\rho(g)^* +
I\right] }\, dg\, d \xi\, d\lambda,
\end{equation*}
 where
$$Y_p=\{\la\in\R^m: \la_1>\cdots >\la_p>0>\la_{p+1}>\cdots
>\la_m\}\;.$$\end{theo}

\subsection{Proof of Lemma \ref{IZint}}\label{SIMPLIFY}

We write
\begin{equation}\label{ical}\ical(z_0)=\frac 1
{\pi^{d_m}\det\Lambda(z_0)} \int_{\sym(m,\C) \times \C}
\left|\det(HH^*-|x|^2I) \right| \exp\left( -{\langle
\La(z_0)^{-1}(H, x),(H, x) \rangle}\right) dH dx\;.\end{equation}
Here, $H$ (previously denoted by $H'$) is a complex $m \times m$
symmetric matrix, so $H^* = \overline{H}$.  The proof is
basically to rewrite \eqref{ical} using the Itzykson-Zuber
integral and Gaussian integration.

We first observe that $$\ical(z_0) = \lim_{\epsilon' \to 0}\;
\lim_{\epsilon \to 0}
\ical_{\epsilon,\ep'} (z_0)\;,$$ where $\ical_{\epsilon,\ep'}(z_0)$ is
the absolutely convergent integral,
\begin{eqnarray}
\ical_{\epsilon,\ep'}(z_0)  & = &
\frac{1}{(2\pi)^{m^2}\pi^{d_m}\det\La}
\int_{{\mathcal H}_m} \int_{\hcal_m} \int_{\sym(m,\C) \times \C}
\left|\det P\right|  e^{- \epsilon  Tr \Xi^*
\Xi-\ep' Tr P^*P } e^{i \langle \Xi, P - HH^* +|x|^2I\rangle}\nonumber\\
&&\quad\times\ \exp\left( -{\langle \La^{-1}(H, x),(H, x)
\rangle}\right)\, dH \,dx \,dP \,d \Xi .\label{ical1a}
\end{eqnarray}
 Here, ${\mathcal H}_m$ is the space of $m \times m$
Hermitian matrices.  Absolute convergence is guaranteed
by the Gaussian factors in each variable $(H, x, \Xi, P)$. If the
$d \Xi$ integral is done first, we obtain a dual Gaussian which
converges (in the sense of tempered distributions) to the delta
function $\de_{HH^* -\half |x|^2}(P)$ as $\ep\to 0$. Then, as
$\epsilon' \to 0$, the $d P$ integral then evaluates the integrand
at $P = H H^*-|x|^2I$ and we retrieve the original integral
$\ical(z_0)$.

We next  conjugate $P$ in (\ref{ical1a}) to a diagonal matrix
$D(\lambda)$ with $\lambda = (\lambda_1, \dots, \lambda_m)$ by an
element $h\in \U(m)$.  Recalling
that
\begin{equation} \label{IDENTITY} \int_{\hcal_m}\phi(P)\,dP =
c_m' \int_{\R^m}\int_{\U(m)}\phi(
hD(\lambda)h^*)\De(\lambda)^2\,dh\,d\la\,, \qquad c_m'= \frac{(2\pi)^{{m\choose 2}}}
{\prod_{j=1}^m j!}  \end{equation}
(see for example \cite[(1.9)]{ZZ}), we then obtain
\begin{eqnarray*}
\ical_{\epsilon,\ep'} & = &  \frac{c'_m}{(2\pi)^{m^2}\,\pi^{d_m}\det\La}
\int_{\U(m) }
\int_{\sym(m,\C)\times\C}
 \int_{{\mathcal H}_m} \int_{\R^m}
 \left|\det(D(\lambda) )  \right| \\ & & \times\
 e^{- \epsilon (Tr D(\lambda)^* D(\lambda) + Tr \Xi^* \Xi)} e^{i \langle \Xi, h D(\lambda) h^* + |x|^2 I  - H^*H
\rangle} \Delta(\lambda)^2\\ & &\times\ \exp\left( -{\langle
\La^{-1}(H, x),(H, x) \rangle}\right)\, d\la\,d\Xi\, dH\, dx \, dh
\end{eqnarray*} Since  the factor $\int_{\U(m) }
e^{i \langle \Xi, h D(\lambda) h^* \rangle} dh$ is invariant under
the conjugation $\Xi \to g^* \Xi g$ with $g \in \U(m)$, we apply
the same identity \eqref{IDENTITY} in the  $\Xi$ variable. We
write $\Xi = g^{-1} D(\xi) g$ where  $D(\xi)$ is diagonal. This
replaces $d \Xi$ by $\Delta(\xi)^2 d \xi$. The inner product is
bi-invariant so we may transfer the conjugation to $HH^*$.  We
thus obtain:
\begin{eqnarray}
\ical_{\epsilon,\ep'} & = &
\frac{(c'_m)^2}{(2\pi)^{m^2}\,\pi^{d_m}\det\La}
\int_{\U(m) }
\int_{\U(m)} \int_{\sym(m,\C)\times\C}
 \int_{\R^m} \int_{\R^m}
 \left|\det(D(\lambda) )  \right| e^{- \epsilon (|\xi|^2 + |\lambda|^2)} \nonumber \\ & & \times\
 \exp\big[i \langle D(\xi),
h D(\lambda) h^* + |x|^2 I  - g HH^* g^{*} \rangle -{\langle
\La^{-1}(H, x),(H, x) \rangle}\big] \nonumber \\ & & \times\
\Delta(\lambda)^2 \Delta(\xi)^2\,d\xi\, d\la\, dH\, dx \, dg
\,dh\;.\label{ical2}
\end{eqnarray}

Next we recognize the integral  $\int_{\U(m)} e^{i \langle D(\xi),
h D(\lambda) h^* \rangle}dh$ as   the well-known
Itzykson-Zuber-Harish-Chandra integral \cite{Ha} (cf., \cite{ZZ}):
\begin{equation}\label{IZ}J(D(\lambda), D(\xi))   =
(-i)^{m(m-1)/2}\left({\textstyle \prod_{j=1}^{m-1}j!}\right)\frac{\det
[e^{i\lambda_j
\xi_k}]_{j,k}}{\Delta(\lambda) \Delta(\xi)}\;.\end{equation}
We note
that both numerator and denominator are anti-symmetric in $\xi_j$
and $\lambda_j$ under permutation, so that the ratio is
well-defined.

We substitute \eqref{IZ} into \eqref{ical2} and expand
$$\det [e^{i
\xi_j \lambda_k}]_{jk} = \sum_{\sigma \in S_m} (-1)^\sigma\;
e^{i\langle  \xi, \sigma(\lambda)\rangle }, $$  obtaining a sum of
$m!$ integrals.  However, by making the change of variables
$\la'=\sigma(\la)$ and noting that $\De(\sigma(\la))=  (-1)^\sigma
\De(\la)$, we see that these integrals are equal, and
\eqref{ical2} then becomes
\begin{eqnarray}
\ical_{\epsilon,\ep'} & = &
\frac{(-i)^{m(m-1)/2}}{(2\pi)^m(\prod_{j=1}^mj!)\pi^{d_m}\det\Lambda}
\int_{\U(m)}
\int_{\sym(m,\C)\times\C}
 \int_{\R^m} \int_{\R^m}\Delta(\lambda) \Delta(\xi)\,
 \left|\det(D(\lambda) )\right|\,
 \nonumber \\ &  \times & e^{i\langle \la,\xi\rangle}  e^{- \epsilon
(|\xi|^2 + |\lambda|^2)}
 e^{i \langle D(\xi),
 |x|^2 I  - g HH^* g^{*} \rangle -{\langle
\La^{-1}(H, x),(H, x) \rangle}} \,d\xi\, d\la\, dH\, dx \,
dg\;. \label{ical3}
\end{eqnarray} Further we observe that the $dH dx$ integral
is a Gaussian integral.  We simplify the phase by noting that
$$\langle D(\xi), gHH^*g^*-|x|^2 I\rangle =
Tr (D(\xi) gHg^t \bar gH^*g^*) - Tr D(\xi)\,|x|^2 = \left\langle
\wh D(\xi)\rho(g)(H,x), \rho(g)(H,x)\right\rangle$$ where $\wh
D(\xi)$ and $\rho(g)$ are as in the statement of the theorem.
Thus,
\begin{eqnarray}\frac 1{\pi^{d_m}\det\La}
\int_{\sym(m,\C)\times\C}
 \exp\big[i \langle D(\xi),
 |x|^2 I  - g HH^* g^{*} \rangle -{\langle
\La^{-1}(H, x),(H, x) \rangle}\big] \,dH\, dx
\hspace{-3in}\nonumber\\& = & \frac
{1}{\det\La\,\det\left[i \rho(g)^* \wh D(\xi)
\rho(g)+\La\inv\right]}\nonumber\\
& = & \frac {1}{\det\left[i \rho(g)^* \wh
D(\xi) \rho(g)\La+I\right]}\nonumber\\ &=&  \frac
{1}{\det\left[i  \wh D(\xi) \rho(g)\La
\rho(g)^* +I\right]}\;. \label{Hxint}\end{eqnarray} Substituting
\eqref{Hxint} into \eqref{ical3}, we obtain
 the desired formula.\qed

\subsection{Proof of Theorem \ref{Iq}}\label{Iq-proof}

The proof is similar to that of Theorem \ref{I} above,
with the modifications provided by Theorem \ref{Morse}.   The proof
of Lemma \ref{IZint} shows that
\begin{eqnarray*}&& \hspace{-.2in}\frac{1}{\det\La}
\int_{{\bf S}_{m,k}} \left|\det(HH{}^*-|x|^2I)\right|e^{
-{\langle \La^{-1}(H, x),(H, x) \rangle}}\,dH\,dx\\&&
=\ \frac{(-i)^{m(m-1)/2}}{(2\pi)^m\prod_{j=1}^mj!} \,\lim_{\ep' \to
0^+}
\int_{Y'_{m-k} } \lim_{\ep \to
0^+} \int_{\R^m} \int_{\U(m)} \!\! \frac{\Delta(\xi)\,
\Delta(\lambda)\; |\prod_j \la_j| \,e^{i\langle \xi, \lambda
\rangle}  e^{- \epsilon |\xi|^2 -\epsilon' |\lambda|^2}}{\det\left[
i\wh D(\xi)\rho(g)\La\rho(g)^* + I\right] }\, dg\, d \xi\, d\lambda,
\end{eqnarray*} where $Y'_p$ denotes the set of points in $\R^m$ with
exactly $p$ coordinates positive. Since the integrand on the right
is invariant under identical simultaneous permutations of the
$\xi_j$ and the $\lambda_j$, it follows that the integral  equals
$m!$ times the corresponding integral over $Y_{m-k}$.  The desired
formula then follows from Theorem \ref{Morse}.\qed

\section{Exact formula for $\CP^m$ }\label{s-exact} To illustrate
our results for fixed $N$,
we compute the density $\kcal^\crit_{N,q}(z)$ of the expected
distribution of critical points of Morse index $q$ of $\log \|s_N\|_{h^N}$
  for random sections $s_N\in H^0(\CP^m, \ocal(N))$, where $h^N$ is the
Fubini-Study metric on $\ocal(N)$.  Here, the probability measure
on $H^0(\CP^m,\ocal(N))$ is  the Gaussian measure induced from
$h^N$ and the volume form $V=\frac 1{m!}\om^m_\FS$ on $\CP^m$.
Since this Hermitian metric and  Gaussian measure are invariant
under the $\SU(m+1)$ action on $\CP^m$, the density is independent
of the point $z\in \CP^m$, and hence the expected number of
critical points of Morse index $q$ is given by
\begin{equation}\label{NK}
\ncal^\crit_{N,q}(\CP^m)=\frac{\pi^m }{m!}
\kcal^\crit_{N,q}(z)\;.\end{equation} These numbers turn out to be
rational functions of $N$ given by the following integral formula:

\begin{prop} \label{numcrit}   The expected number of critical points  of Morse
index $q$ for random sections $s_N\in H^0(\CP^m, \ocal(N))$ is
given by   \begin{eqnarray*}
{\mathcal  N}^\crit_{N,q}(\CP^m) &=& \frac{ 2^{\frac
{m^2+m+2}2}}{\prod_{j=1}^m j!}\ \frac{(N-1)^{m+1}
}{(m+2)N-2}
\\ & &\times
\int_{Y_{2m-q} }d\la\,
\left|{\textstyle \prod_{j=1}^m\lambda_j} \right| \,\Delta(\lambda)\,e^{
-\sum_{j=1}^m \la_j} \times
\left\{\begin{array}{ll}\displaystyle \, e^{(m+2-2/N) \lambda_m}
&\mbox{for
}\ q>m\\[10pt]
\displaystyle 1 &\mbox{for }\ q=m\end{array}\right\},
\end{eqnarray*} for $N\ge 2$, where
$Y_{2m-q}$ is as in Theorem \ref{Iq}.
\end{prop}

Thus, the expected number ${\mathcal  N}^\crit_{N,m}(\CP^m)$ of critical
points of minimum Morse index is of the form $\kappa_m\,\frac{(N-1)^{m+1}
}{(m+2)N-2}$, for some constant $\kappa_m>0$.  It is also easy to see from
the form of the above integral that  the
expected number ${\mathcal  N}^\crit_{N,q}(\CP^m)$ of critical points of
each Morse index $q$ is a rational function of $N$.  In Appendix~1, we give
explicit formulas for ${\mathcal  N}^\crit_{N,q}(\CP^m)$ in low dimensions
obtained by using Maple to evaluate the integral in Proposition
\ref{numcrit}.

\subsection{Itzykson-Zuber formula on $\CP^m$}

The following formula derived from Theorem \ref{Iq} is the
starting point for our proof of Proposition \ref{numcrit}.

\begin{lem} \label{kcritPm}  The expected critical point density for random
sections $s_N\in H^0(\CP^m, \ocal(N))$ of Morse index $q$
is given by   \begin{eqnarray*}
\kcal^\crit_{N,q}(z) &=& i^{m+1}\frac{m!\,|c_m|}{N^m}\lim_{\ep'
\to 0^+} \int_{Y_{2m-q} }d\la\, \left|{\textstyle
\prod_j\lambda_j} \right| \,\Delta(\lambda)\,e^{ -\ep'
|\lambda|^2}\\&&\times\ \lim_{\ep \to 0^+} \int_{\R^m}\frac {
 \Delta(\xi)\,e^{i\langle \lambda,\xi\rangle}e^{-
\epsilon |\xi|^2}\,d\xi}{\left(N^2\sum\xi_j+i\right)\prod_{1\le j \le k\le m}
\{i-N(N-1)(\xi_j+\xi_k)\}}\,,\end{eqnarray*} where $c_m$ and $Y_{2m-q}$
are as in Theorems \ref{I} and \ref{Iq}.
\end{lem}

\begin{proof}
Since the critical point  density $\kcal^\crit_{N,q}$ is constant,
it suffices to compute it  at $z=0\in\C^m\subset \CP^m$, using the
local frame $e_L$ corresponding to the homogeneous (linear)
polynomial $z_0$. We recall that the \szego kernel is given by
$$ \Pi_{H^0(\CP^m,\ocal(N))}(z,w) = \frac {(N+m)!}{\pi^m N!}
(1+z\cdot \bar w)^N e_L(z)\otimes \overline{e_L(w)}\;.$$ (See, for
example, \cite[\S 1.3]{BSZ1}.) Since the formula in Theorem
\ref{KNcrit2} is invariant when the \szego kernel is multiplied by a
constant, we can replace the above by the {\it normalized \szego
kernel}
\begin{equation}
\label{NSK} \wt \Pi_N(z,w):=(1+z\cdot \bar w)^N \end{equation} in
our computation.

We notice that
$$K(z)=-\log\|e_L(z)\|^2_h=\log (1+\|z\|^2)\;,$$
$$K(0) = \frac {\d K}{\d z} (0)
=\frac {\d^2K}{\d^2 z} (0) = 0\;. $$ Hence when computing the
(normalized) matrices $\wt A_N,\ \wt  B_N,\ \wt C_N$, we can take the
usual derivatives of $\wt\Pi_N$.
Indeed, we have
\begin{eqnarray*} \frac{\d\wt\Pi_N}{\d{z_j}} &=& N(1+z\cdot \bar w)^{N-1}
\bar w_j\;,\\
\frac{\d^2\wt\Pi_N}{\d{z_j}\d\bar w_{j'}} &=&\de_{jj'}
N(1+z\cdot\bar w)^{N-1} +
N(N-1)(1+z\cdot\bar w)^{N-2}z_{j'}\bar w_j\;,\\
\frac{\d^4\wt\Pi_N}{\d{z_j}\d z_q\d\bar w_{j'}\d\bar w_{q'}}
(0,0)&=&N(N-1)(\de_{jj'} \de_{qq'} + \de_{j'q}\de_{jq'})\;.
\end{eqnarray*} It follows that
\begin{equation}\label{ABC} \wt A_N=NI\;,\quad \wt B_N= 0\;,\quad \wt\La_N=\wt C_N=
\begin{pmatrix} 2N(N-1)\hat I
&0\\ 0 & N^2 \end{pmatrix}\;,
\end{equation} where $\hat I$ is the identity matrix of rank
${m+1\choose 2}$.

The stated formula now follows from Theorem \ref{Iq} by observing that
$\rho(g)\wt\La_N\,\rho(g)^*=\wt\La_N$, and
$$ \det\left[ i\wh
D(\xi)\wt\La_N + I\right]=(- i)^{\frac{m^2+m+2}2}
\left(N^2\sum\xi_j+i\right)\prod_{1\le j \le k\le m}
\{i-N(N-1)(\xi_j+\xi_k)\}\;.$$
\end{proof}

\subsection{Evaluating the inner integral by residues} \label{residues}

To complete the proof of Proposition~\ref{numcrit}, we must evaluate the
$d
\xi$ integral in Lemma
\ref{kcritPm}. We begin by writing
\begin{equation}\label{kcrit*}\kcal^\crit_{N,q}(z) =
i^{m+1}{m!\,|c_m|}\lim_{\ep'
\to 0^+} \int_{Y_{2m-q} }d\la\, \left|{\textstyle
\prod_j\lambda_j} \right| \,\Delta(\lambda)\,e^{ -\ep'
|\lambda|^2} \ical_{N,\lambda}\;,\end{equation} where
$$\ical_{N,\la}:=\lim_{\ep \to 0^+} \frac 1{N^m}\int_{\R^m}\frac {
 \Delta(\xi)\,e^{i\langle \lambda,\xi\rangle}e^{-
\epsilon |\xi|^2}\,d\xi}{\left(N^2\sum\xi_j+i\right)\prod_{1\le j \le k\le m}
\{i-N(N-1)(\xi_j+\xi_k)\}}.$$  To simplify the constant factors, we make the
redefinitions $\xi_j = (t_j+i)/2N(N-1)$ and
$\lambda_j\rightarrow 2N(N-1)\lambda_j$, after which \eqref{kcrit*} holds
with
$$
\ical_{N,\la} =(-1)^{\frac {m(m+1)}2}\, 2^{\frac
{(m+1)(m+2)}2}\,\frac{(N-1)^{m+1} }{N}
e^{-\sum \la_j}\,\ical(\la;c)\;,$$
where
\begin{eqnarray}\label{Ilambda} \ical(\la;c) &=&\lim_{\ep\to 0^+}
\int_{(\R-i)^m}\frac {
\Delta(t)\,e^{i\langle \lambda,t\rangle}\,e^{-\ep \sum |t_j|^2}}
 {\left(\sum t_j+ic\right)\prod_{1\le j \le
k\le m} (t_j+t_k)}\,dt\;,\\
 c &=& m+2-2/N\;.\label{c}\end{eqnarray}
Recalling the definition of $c_m$ in the statement of Theorem \ref{I}, we
therefore have:
\begin{eqnarray}\label{kcrit}\kcal^\crit_{N,q}(z) &=&
i^{(m+1)^2}\, \frac{ 2^{\frac
{m^2+m+2}2}}{\pi^{2m}\,\prod_{j=1}^{m-1}j!}\,\frac{(N-1)^{m+1}
}{N}\nonumber \\ &&\times\lim_{\ep'
\to 0^+} \int_{Y_{2m-q} }d\la\, \left|{\textstyle
\prod_j\lambda_j} \right| \,\Delta(\lambda)\,e^{ -\ep'
|\lambda|^2}
e^{-\sum \la_j}\,\ical(\la;c)\;,\end{eqnarray} where $\ical(\la;c)$ is given
by \eqref{Ilambda}--\eqref{c}. Proposition \ref{numcrit} now follows from
\eqref{NK},
\eqref{kcrit}, and the following lemma with $p=2m-q$, $c=m+2-2/N$.

\begin{lem} \label{xiint} Let $0\le p\le m$ and let $c>0$.  Then for
$$\la_1> \cdots >\lambda_p
> 0 > \lambda_{p+1}>\cdots
>\la_m\;,$$
we have
$$
 \ical(\la;c)=\left\{\begin{array}{ll}\displaystyle
i^{m^2-1}\,\frac{\pi^m}{c}\, e^{c\lambda_m} \quad &\mbox{for
}\ p<m\\[10pt]
\displaystyle i^{m^2-1}\,\frac{\pi^m}{c} &\mbox{for }\ p=m
\end{array}\right.\ ,$$ where $\ical(\la,c)$ is given by \eqref{Ilambda}.
\end{lem}
\noindent  We note that $\ical(\la,c)$, resp.\ $i\,\ical(\la,c)$, is positive
if $m$ is  odd, resp.\ even, and thus by monotone
convergence we can set $\ep'=0$ in \eqref{kcrit}.
\begin{proof} We let

\begin{equation}\label{integrand}
{\mathcal I}(\lambda,t;c) =
\frac {\Delta(t)\,e^{i\langle \lambda,t\rangle}}
 {\left(\sum t_j+ic\right)\prod_{1\le j \le
k\le m} (t_j+t_k)}\,dt\;,
\qquad \mbox{for }\ c>0\;. \end{equation}

We note that $\int_{(\R-i)^m}{\mathcal
I}(\lambda,t;c)\,dt$ is a
tempered distribution (in $\la$). Furthermore, the map
$$(\ep_1,\dots,\ep_m)\mapsto
\int_{(\R-i)^m}{\mathcal
I}(\lambda,t;c)\,e^{ -\sum\ep_j|t_j|^2}\,dt$$ is a continuous map from
$[0,+\infty)^m$ to the tempered distributions. Here,
$(\R-i)^m$ is the change in the contour $\R^m$  obtained by
translating $\R \to \R - i$ in each factor.  Hence
\begin{equation}\label{it}\ical(\la;c)=\int_{(\R-i)^m}{\mathcal
I}(\lambda,t;c)\,dt = \lim_{\ep_m\to 0^+} \cdots \lim_{\ep_1\to
0^+} \int_{(\R-i)^m}{\mathcal I}(\lambda,t;c)\,e^{
-\sum\ep_j|t_j|^2}\,dt\;.\end{equation}

We now
use \eqref{it} to evaluate ${\mathcal I}(\lambda;c)$ by iterated
residues.
We first suppose that $p>0$, and we start by doing the
integral over
$t_1$.  Since the $t_1$ integral is absolutely convergent when $\ep_1=0$, we
can set
$\ep_1=0$ and do the integral by residues. If
$p>0$ we close the contour in the upper half plane, and pick up poles at
$t_1=0$, and at
$t_1=-t_j$ for
$j\ne 1$.  The pole at
$t_1=-ic-\sum_{j\ne i}t_j$ is below the contour.

The residue of $\ical(\la,t;c)$ at the pole $t_1=0$
is
\begin{equation}\label{res1}\frac{(-1)^{m-1}}{2} {\mathcal
I}(\lambda_2,\ldots,\la_m,t_2,\dots,t_m;c)\;.\end{equation}
The residue at the pole $t_1=-t_2$ is
\begin{eqnarray*}&&\hspace{-.3in}
\frac{\pm e^{i[(\lambda_2-\lambda_1)t_2+\la_3t_3+\cdots\la_mt_m]}
2t_2(t_2+t_3)\cdots(t_2+t_m)\Delta(t_2,\dots,t_m)} {(t_3+\cdots
+t_m+ci)\,2t_2(-t_2+t_3)\cdots(-t_2+t_m)\prod_{2\le j \le
 k\le m} (t_j+t_k)}\\ &&=\ \frac{\pm
e^{i(\lambda_2-\lambda_1)t_2}\,e^{-\ep_2|t_2|^2}}{2t_2}\,
\ical(\la_3,\dots,\la_m,t_3,\dots,t_m;c)\;.
\end{eqnarray*}
When we then do the $t_2$ integral  and let $\ep_2\to 0^+$, we get zero. Indeed,
$$\lim_{\ep\to 0+}\int_{\R-i}\frac{
e^{i(\lambda_2-\lambda_1)t_2}\,e^{-\ep_2|t_2|^2}}{2t_2}\,dt_2=0\;,$$
since $\la_2-\la_1<0$ and the pole at $t_2=0$ is above the contour.
Similarly, when we compute the residue of the pole
$t_1=-t_j$,
$j>2$, and then perform the $t_j$ integration, we also get zero.  Hence we
can ignore the residues of the poles $t_1=-t_j$.

Applying \eqref{res1} recursively, the integral with $p>0$ can be reduced to
the case with all $\la$'s negative:
\begin{equation}\label{p=0}
{\mathcal I}(\lambda;c) =
(-1)^{(m-1)+(m-2)+\cdots+(m-p)}(\pi i)^p\, {\mathcal
I}(\lambda_{p+1},\ldots,\lambda_m;c)\; .
\end{equation}

We now treat the case $p=0$ (i.e., $0>\la_1>\cdots>\la_m$). This
time, we do the $t_m$ contour integral first. We close it in the
lower half plane, picking up the residue at $t_m=-ic-\sum_{1\le
k<m}t_k$. This residue is
\begin{eqnarray}&&\hspace{-.3in}
\rcal(\la_1,\dots,\la_{m-1},t_1,\dots,t_{m-1};c): =\nonumber \\ &&
\frac{\Delta(t_1,\ldots,t_{m-1}) \prod_{k< m} (ic+\sum_{l<m} t_l
+t_k) e^{c\lambda_m + i\sum_j (\lambda_j-\lambda_m)t_j}}
{2(-ic-\sum_{l<m} t_l) \prod_{1\le j\le k\le m-1} (t_j+t_k)
\prod_{k<m} (-ic-\sum_{l<m,l\ne k} t_l) }\;.\label{res2}
\end{eqnarray} (To simplify the
discussion, we set
$\ep=0$, and regard the integrals as distributions, as above.)
Next we perform the $t_1$ integration. Since $\lambda_m$ is the
most negative eigenvalue, we close the contour  in the upper half
plane. The terms in the denominator with $ic$ all give poles in
the lower half plane, so can be ignored.  And, the poles
$t_1=-t_j$ will be ignorable, by the same type of reasoning we saw
earlier. Indeed, after computing the residue at $t_1=-t_j$ we find
that $t_j$ appears in the exponent as $e^{i(\la_j-\la_1)t_j}$ with
$\la_j-\la_1<0$ and the only factor of the denominator with a zero
below the contour is $ic+t_2+\cdots +t_m$; but this factor also
appears in the numerator and hence the $t_j$ integral gives zero.

This leaves the residue at all $t_j=0$ with $1\le j \le m-1$.
 The residue at $t_1=0$ of   (\ref{res2}) is

$$ \mbox{Res}\;|_{t_1 =0}\; \rcal(\la_1,\dots,\la_{m-1},t_1,\dots,t_{m-1};c) = \frac {(-1)^{m-1}}{2}
\rcal(\la_2,\dots,\la_{m-1},t_2,\dots,t_{m-1};c)\;.$$
Continuing recursively, for the case $p=0$, we obtain (remembering
that the $t_m$ pole below the contour contributes negatively):
\begin{equation}\label{Ip=0} \ical(\la;c)=
(-1)^{m(m-1)/2}\,(\pi i)^m\,\left(\frac{-i}{c}\right)\,
e^{c\lambda_m}\;.\end{equation}

Combining \eqref{Ip=0} (with $m$ replaced by $m-p$) and \eqref{p=0},
we obtain the formula of the proposition.\end{proof}

\subsection{Dimensional dependence: a conjecture}

In this article, we are studying the  $N$-dependence of the
critical point density and number of critical points, but  the
growth rate of  ${\mathcal N}^\crit_{N, q}$  (and the
growth rate of the total number of critical points) as the
dimension $m \to \infty$ is of considerable interest as well.  In
the vacuum statistics problem in string/M theory, one is
interested in the total number of critical points of certain
holomorphic sections (known as flux superpotentials) of a line
bundle $\lcal \to \ccal$ over a moduli space $\ccal$ of complex
structures on a Calabi-Yau fourfold $X \times T^2$ where $X$ is a
Calabi-Yau threefold and $T^2$ is an elliptic curve. The relevant
dimension is $m = b_3/2 -1$ where $b_3$  is  the third Betti
number of $X$.  As discussed in some detail in \S 7.3 of
\cite{DSZ2}, vacuum statistics involves integrals like (\ref{b0})
over certain subspaces of $\sym(b_3/2 -1, \C)$. The growth rate in $b_3$
of the number of vacua is important to obtain an order of
magnitude of string vacua,  since $b_3$ is rather large for
typical Calabi-Yau threefolds. In \cite[\S 7.3]{DSZ2}, we state a
conjectural formula for the growth in $b_3$ of such integrals
which is based on the explicit calculations in this article for
$\ocal(N) \to \CP^m$ for fixed $N$.

We recall that by Proposition \ref{numcrit}, the expected number of critical
points on $\CP^m$ of minimum Morse index is of the form
$\kappa_m\,\frac{(N-1)^{m+1}
}{(m+2)N-2}$. Our Maple computations given in
Appendix~1 show that $\kappa_m=2(m+1)$, and  furthermore that the leading
terms of the expansion are monotonically decreasing in the Morse index $q$,
for each
$m\le 6$; in particular, on average there are more critical points on $\CP^m$
of Morse index $m$ than there are of any index $q>m$ (at least for
$m\le 6$). Recall that by our universality results (Theorems
\ref{UNIVDIST}--\ref{univmorse} and Corollary \ref{EN}), the leading
coefficient in the
$N$-expansion of
$\ncal_{N,q}(\CP^m)$ is the universal coefficient
$\frac{\pi^m}{m!}\,b_{0q}(m)$.  Thus we make the following conjecture
based on our Maple computations:

\begin{conj} Let $n_q(m):=\frac{\pi^m}{m!}\,b_{0q}(m)$ denote the
leading coefficient in the expansion of $\ncal^\crit_{N,q,h}$ from
Corollary
\ref{EN}, and let $n(m)=\sum_{q=m}^{2m} n_q(m)=
\frac{\pi^m}{m!}\,b_{0}(m)$: $$\ncal^\crit_{N,q,h}
\sim n_q(m)\,c_1(L)^m\,N^m ,\qquad \ncal^\crit_{N,h}
\sim n(m)\,c_1(L)^m\,N^m .$$  Then $$n_m(m)=
2\,\frac{m+1}{m+2} > n_{m+1}(m) >
\cdots > n_{2m}(m)\;,$$ and hence the expected total number of critical
points
$$\ncal^\crit_{N,h} \sim
n(m)\,c_1(L)^m\,N^m \;,\qquad 2\,\frac{m+1}{m+2} <n(m) <
2\,\frac{(m+1)^2}{m+2}\;.$$
\end{conj}

Equivalently, $\ncal^\crit_{N,h} \sim n(m)\,\deg L^N$, where $\deg
L^N$ is the number of simultaneous zeros of $m$ generic sections of
$L^N$.  This conjecture implies that the expected number of critical
points of Morse index $m$ grows exponentially in the dimension, a
growth rate consistent with quite analogous estimates of vacua in
string/M theory
\cite{DSZ2, DD} and of metastable states of spin glasses \cite{Fy}.

\section{Asymptotics of the expected number of critical
points}

In this section, we compute the asymptotics of the expected
density and number of critical points of sections of powers $L^N$
of a positive holomorphic line bundle. In particular, we prove
 Theorems \ref{UNIVDIST},  \ref{critRS}, and \ref{univmorse} as well as Corollary \ref{EN}.

 \subsection{Proof of Theorems \ref{UNIVDIST}--\ref{univmorse}}

 We begin with some further background on the \szego kernel.

\subsubsection{\szego kernel asymptotics}

We first use the asymptotic expansion of the \szego kernel to show
that $\kcal^\crit_N(z)$ has an expansion of the type given in
Theorems~\ref{UNIVDIST}--\ref{univmorse}.
 It is evident from
Theorem \ref{KNcrit2} and Theorem \ref {Morse}, respectively, and
from formulas (\ref{LambdaN})--(\ref{CN}) for $A_N$ and $\Lambda_N$
that the  asymptotics of the critical-point densities
$\kcal_N^{\crit}(z)$ and $\kcal_{N,q}^{\crit}(z)$, respectively,
can be determined by canonical algebraic operations on the
asymptotics of the following derivatives of the \szego kernel (
$1\le j\le m\,, 1\le j\le q\le m\,, 1\le j'\le q'\le m$)
\begin{itemize}

\item $\nabla_{z_j}\Pi_N(z, w)|_{z = w}$;

\item $\nabla_{z_j}\nabla_{\bar w_{j'}}\Pi_N(z,w)|_{z = w};$

\item $ \nabla_{z_q}\nabla_{z_j}\Pi_N(z, w)|_{z = w}$ and
 $\nabla_{\bar w_{q'}}\nabla_{\bar w_{j'}}\Pi_N(z, w)|_{z =
w};$

\item
$ \nabla_{z_j}\nabla_{\bar w_{q'}}\nabla_{\bar w_{j'} }\Pi_N(z,
w)|_{z = w}$;

\item $\nabla_{z_q}\nabla_{z_j}\nabla_{\bar w_{q'}}\nabla_{\bar w_{j'}}\Pi_N(z, w)|_{z = w}$

\end{itemize}

(Here we write $\Pi(z,w)=\Pi(z,0;w,0)$.) We can obtain their
asymptotics by differentiating the following Tian-Yau-Zelditch
expansion:

\begin{theo} \cite{Ze, T1, Y2} \label{tyz}
Let  $(L,h)\rightarrow M$ be a positive Hermitian holomorphic
line bundle over a compact complex manifold $M$  of dimension $m$ with \kahler form
$\om_h=\frac i2 \Theta_h$. Then there is a complete asymptotic
expansion:
\begin{equation}\label{fud1}
\Pi_N(z, z) \sim \frac{N^m}{\pi^m}\left[1+a_1(z) N^{-1}
+a_2(z) N^{-2}+\cdots\right]\,,
\end{equation}
for certain smooth coefficients $a_j(z)$.
\end{theo}

To apply \eqref{fud1} to the differentiated \szego kernel, we use
\eqref{nablaz}--\eqref{nablaw}. By a change of frame in $L$, we can assume that $K$ and
its holomorphic derivatives up to any fixed order, as well as the anti-holomorphic
derivatives,
 vanish at $z_0$. Writing $\d_j=\d/\d z_j$, we then  have:
\begin{eqnarray} &&\nabla_{z_j}\Pi_N(z_0,z_0) = \frac{\d F_N}{\d
z_j}(z_0,\bar z_0) = \d_jF_N(z,\bar z)|_{z_0} =
\d_j\left.\left[e^{NK(z)}\Pi_N(z,z)\right]\right|_{z_0}\,,\nonumber\\
&&\nabla_{z_j}\nabla_{\bar w_{j'}}\Pi_N(z_0,z_0) = \frac{\d^2F_N}{\d
z_j\d\bar w_{j'}}(z_0,\bar z_0) = \d_j\bar\d_{j'}F_N(z,\bar z)|_{z_0} =
\d_j\bar\d_{j'}
\left.\left[e^{NK(z)}\Pi_N(z,z)\right] \right|_{z_0}\,,\nonumber\\[-5pt]
&& \qquad\vdots\label{derivatives}\\
&&\nabla_{z_j}\nabla_{z_q}\nabla_{\bar w_{j'}}\nabla_{\bar
w_{q'}}\Pi_N(z_0,z_0)  = \d_j\d_q\bar\d_{j'} \bar\d_{q'}
\left.\left[e^{NK(z)}\Pi_N(z,z)\right]
\right|_{z_0}\,.\nonumber\end{eqnarray}  Here, we used the fact
that
$F_N(z,\bar w)$ is holomorphic in $z$ and anti-holomorphic in $w$. (In these expressions,
we have no $\nabla_{\bar{z}_k}$ or
$\nabla_{w_j}$ derivatives of $\Pi_N(z,w)$.)

It follows by substituting \eqref{fud1} into \eqref{derivatives}
that the components of $A$ and $\Lambda$ have asymptotic
expansions in powers of $N$, and hence by Theorem~\ref{KNcrit2},
resp.\ Theorem \ref{Morse}, that $\kcal^\crit_N(z)$, resp.\
$\kcal^\crit_{N,q}(z)$, does. Next we study the coefficients
$b_0,b_1,b_2$ of the expansion of $\kcal^\crit_N(z)$.

\subsubsection{The first three terms of the
expansion}\label{s-expansion}

Integrating the density of critical points, we find that the
expected total number of critical points has the expansion
$$N^{-m}\;\ncal^\crit_{N,h} = \frac{\pi^m}{m!}\,b_0\,
c_1(L)^m +  N^{ - 1}\int_M b_1 dV_h
+ N^{-2} \int_M b_2 dV_h + O(N^{ - 3})\,.$$ The leading
order term is universal.

 We will use Theorem \ref{KNcrit2} and the following result
of Z. Lu \cite{Lu} to calculate the coefficients in these expansions:

\begin{theo}\label{zlu} \cite{Lu}
With the notation as in Theorem \ref{tyz}, each coefficient
$a_j(z)$ is a polynomial of the curvature
 and its covariant derivatives  at $x$.  In particular,
\[
\left\{
\begin{array}{l}
a_1=\frac 12 \rho\\
a_2=\frac 13\Delta\rho+\frac{1}{24}(|R|^2-4|Ric|^2+3\rho^2)\\
\end{array}
\right.
\]
where $R, Ric$ and $\rho$ denotes the curvature tensor, the
Ricci curvature and the scalar curvature of $\om_h$, respectively, and
$\Delta$ denotes the Laplace operator of $(M,\om_h)$.
\end{theo}

We now calculate $A_N$ and $\Lambda_N$ to two orders. The key
point is to calculate the mixed derivatives of $\Pi_N$ on the
diagonal. It is convenient to do the calculation in K\"ahler
normal coordinates about a point $z_0$ in $M$.

It is well known that in terms of \kahler normal coordinates $\{z_j\}$,
the K\"ahler potential $K$ has the expansion:
\begin{equation}\label{potexp}K(z,\bar z)= \|z\|^2 -\frac 14 \sum R_{j\bar
kp\bar q}(z_0)z_j\bar z_{\bar k} z_p\bar z_{\bar q} +
O(\|z\|^5)\;.\end{equation} (In general, $K$ contains a
pluriharmonic term $f(z) +\overline{f(z)}$, but  a change of frame for $L$
eliminates that term up to fourth order.)

We further use the notation $K_j = \d_j
K,\ K_{\bar{j}} = \dbar_j K$.
 We first claim that
 \begin{equation} \label{Aexp}A = N I + a_1 I+
N^{-1}\left\{a_2I+\big(\d_j\bar\d_{j'}a_1\big)\right\}+\cdots.
\end{equation}
Indeed,
by \eqref{derivatives},
\begin{eqnarray*}
\d_j\left[e^{NK(z)}\Pi_N(z,z)\right]
&=& e^{NK}
\left[NK_j(1+a_1N\inv+a_2N^{-2}) +\d_ja_1N\inv
+\d_ja_2N^{-2}+\cdots\right]\\
\d_j\dbar_{j'}\left[e^{NK(z)}\Pi_N(z,z)\right]
&=& e^{NK}\big[N^2K_jK_{\bar j'}(1+a_1N\inv+a_2N^{-2}) +K_{\bar
j'}\d_ja_1 +K_{\bar j'}\d_ja_2N^{-1}\\
&& + NK_{j\bar j'}(1+a_1N\inv+a_2N^{-2})
+K_j\dbar_{j'}a_1 +K_{j}\dbar_{j'}a_2N^{-1} \\&&+ \d_j\dbar_{j'} a_1
N\inv\cdots\big].\end{eqnarray*}
Evaluating at $z_0$ using \eqref{potexp}, we then obtain \eqref{Aexp}.

We now compute the expansion of $\Lambda$.  Continuing the above
computation,
\begin{eqnarray*}
\d_j\dbar_{j'}\dbar_{q'}\left[e^{NK(z)}\Pi_N(z,z)\right]
&=& e^{NK}\big[N^2(K_{j\bar q'}K_{\bar j'} + K_{j\bar j'}K_{\bar
q'})(1+a_1N\inv+a_2N^{-2})\\&& +K_{\bar j'}\dbar_{q'}a_1 +K_{\bar
q'}\d_j\dbar_{j'}a_1 +K_{j\bar j'}\dbar_{q'}a_1 +K_{j\bar
q'}\dbar_{j'}a_1 \\ && + NK_{j\bar j' \bar
q'}(1+a_1N\inv)\cdots\big] +\mbox{unimportant terms}.\end{eqnarray*}
(The `unimportant terms' are those which vanish at $z_0$ and whose
holomorphic derivatives also vanish at $z_0$.) We have
\begin{eqnarray}B(z_0) &=&
\Big[ \big(\nabla_{z_j}\nabla_{\bar w_{j'}}\nabla_{\bar
w_{q'}}\Pi_N(z_0,z_0)\big)\quad \big(N\nabla_{z_j}\Pi_N(z_0,z_0)
\big)\Big]\nonumber\\& =& \Big[ \big(\de_{jj'}\dbar_{q'}a_1 +
\de_{jq'}\dbar_{j'}a_1\big)\quad \big(\de_{j}a_1\big) \Big] \ +\
O(N\inv)\;.\label{Bexp}\end{eqnarray}

Differentiating again and evaluating at $z_0$ using \eqref{potexp}, we
obtain
\begin{eqnarray*}
\d_j\dbar_{j'}\dbar_q\dbar_{q'}\left[e^{NK(z)}\Pi_N(z,z)\right]\big|_{z_0}
&=&
\big[N^2(\de_{jj'}\de_{qq'}+\de_{jq'}\de_{qj'})(1+a_1N\inv+a_2N^{-2})\\&&
+\de_{jj'}\d_q\dbar_{q'}a_1+\de_{qq'}\d_j\dbar_{j'}a_1
+\de_{jq'}\d_q\dbar_{j'}a_1 +\de_{qj'}\d_q\dbar_{j'}a_1
\\ && + NK_{j\bar j' q \bar q'}(1+a_1N\inv)\cdots\big]\big|_{z_0}
\;.\end{eqnarray*} Noting that  $K_{j\bar j' q \bar q'}|_{z_0} =
- R_{j\bar j' q \bar
q'} (z_0)$, and recalling that $\La=C-B^*A\inv B$, where
$$C=  \left[
\begin{array}{cc}\big(\nabla_{z_q}\nabla_{z_j}\nabla_{\bar w_{q'}}
\nabla_{\bar w_{j'}}\Pi_N \big)&
\big(N\nabla_{z_q}\nabla_{z_j}\Pi_N\big) \\[8pt]
\big(N\nabla_{\bar w_{q'}}\nabla_{\bar w_{j'}}\Pi_N \big)&
N^2\Pi_N \end{array}\right]\,,$$ we obtain
$$\La(z_0) = N^2\La_0^\half\left(I +N\inv\La_{-1} +N^{-2}\La_{-2}+
\cdots\right)\La_0^\half\;,$$ with
\begin{equation}\label{L0}
\La_0=\begin{pmatrix} 2\hat I  &0\\ 0 &1\end{pmatrix}\;,\end{equation}
\begin{equation}\label{L1}
\La_{-1}=\begin{pmatrix} a_1\hat I - \half\big( R_{j\bar j' q \bar
q'} \big) &0\\ 0
&a_1\end{pmatrix}\;,\end{equation}
\begin{equation}\label{L2}
\La_{-2}=\begin{pmatrix} a_2\hat I + P -\frac{a_1}{2}\big( R_{j\bar j' q
\bar q'} \big)\  &\frac 1{\sqrt 2}\big(\d_j\d_qa_1\big)\\
\frac 1{\sqrt 2} \big(\dbar_j\dbar_qa_1\big)
&a_2\end{pmatrix}\;,\end{equation} where $\hat I$ is the identity operator on
$\sym(m,\C)$, and
$$P=\half\big(\de_{jj'}\d_q\dbar_{q'}a_1 +\de_{qq'}\d_j\dbar_{j'}a_1
+\de_{jq'}\d_q\dbar_{j'}a_1 +\de_{qj'}\d_q\dbar_{j'}a_1\big)\;.$$

To prove Theorem \ref{UNIVDIST}, we want the asymptotics of
$$\kcal^\crit_N (z_0) =
\frac{\pi^{-{m+2\choose 2}}N^m}{\det A_N \det\La_N} \int_{\sym(m,\C)
\times \C} \left|\det(H'H'{}^*-|x|^2I)\right|e^{ -{\langle
\La_N(z_0)^{-1}(H', x),(H', x) \rangle}}\,dH'\,dx\,. $$
Making the change of variables $H'\mapsto \sqrt2 \, N H',\ x\mapsto Nx$,
the integral is transformed to
\begin{equation}\kcal^\crit_N (z_0) =
\frac{\pi^{-{m+2\choose 2}}N^m}{\det \wt A\, \det\wt\La} \int_{\sym(m,\C)
\times \C} \left|\det(2H'H'{}^*-|x|^2I)\right|e^{ -{\langle
\wt\La^{-1}(H', x),(H', x) \rangle}}\,dH'\,dx\,, \end{equation}
where
\begin{equation}\label{newA}\wt A= N^{-2}A_N(z_0), \quad \wt \La =
N^{-2}\La_0^{-\half}\La_N(z_0)\La_0^{-\half}= \left(I +N\inv\La_{-1}
+N^{-2}\La_{-2}+
\cdots\right)\;.
\end{equation}

Next we observe that
$$\wt\Lambda^{-1} = I - \frac{1}{N} \Lambda_{-1} + \frac{1}{N^2} [-
\Lambda_{-2} + \Lambda_{-1}^2 ] + \cdots $$ hence
$$\begin{array}{lll} e^{- \langle \wt\Lambda^{-1} H, H \rangle } &\sim &
e^{-
\langle H, H
\rangle} e^{\left\langle \left [ \frac{1}{N} \Lambda_{-1} + \frac{1}{N^2}
(\Lambda_{-2} - \Lambda_{-1}^2 )\right] H, H \right\rangle} \\ & & \\
& = & e^{- \langle H, H \rangle} \{1 + \frac{1}{N} \langle
\Lambda_{-1} H, H \rangle \\ \\
&+ &\frac{1}{N^2} [  \langle  \Lambda_{-2} H, H \rangle +
\frac{1}{2} \langle \Lambda_{-1} H, H \rangle^2 - \langle
\Lambda^{2}_{-1} H, H \rangle ] \}.
\end{array}$$

Furthermore $$\det \wt\La\inv =  1 -(\tr  \La_{-1})N\inv + \left[\half
\tr (\La_{-1}^2)+\half (\tr  \La_{-1})^2-\tr  \La_{-2}\right]N^{-2} \cdots\,,$$
and similarly for $\det A\inv$.
Altogether, we obtain:
$$\begin{array}{lll} \kcal^\crit_{N} (z) &\sim &
\pi^{-{m+2\choose 2}}N^m\left\{  1 + \frac{1}{N}\left[-  \tr  A_{-1}  -  \tr
\Lambda_{-1}\right]
\right. \\
\\ &+ &\frac{1}{N^2} \left[\half
\tr (\La_{-1}^2)-\left. \tr  \La_{-2} +\half
\tr (A_{-1}^2)-\tr  A_{-2}+\half (\tr  A_{-1}+\tr \La_{-1})^2 \right]\right\}
\\ & &
\\ &\times&  \int_{\sym(m,\C) \times \C}
\left|\det(2H'H'{}^*-|x|^2I)\right|e^{- \langle H, H \rangle} \{1 + \frac{1}{N} \langle
\Lambda_{-1} H, H \rangle \\ \\
&+ &\frac{1}{N^2} [  \langle  \Lambda_{-2} H, H \rangle +
\frac{1}{2} \langle \Lambda_{-1} H, H \rangle^2 - \langle
\Lambda^{2}_{-1} H, H \rangle ] \}\,dH'\,dx\,. \end{array}$$

Expanding, we obtain \begin{equation}
\label{beta}\begin{array}{rcl} \kcal^\crit_{N} (z) &\sim &
b_0N^m+b_1N^{m-1}
+b_2N^{m-2} +\cdots,\\[6pt]
b_0&=&\int d\mu\,,\\[6pt]
b_1 &=& \int\left[\langle \Lambda_{-1} H, H \rangle-  \tr  A_{-1}  -  \tr
\Lambda_{-1}\right]
d\mu\,, \\[6pt] b_2 &= & \int\left[
\half
\tr (\La_{-1}^2)- \tr  \La_{-2} +\half
\tr (A_{-1}^2)-\tr  A_{-2}+\half (\tr  A_{-1}+\tr \La_{-1})^2  \right.\\[6pt]
 && \left. -(\tr
A_{-1}  +  \tr  \Lambda_{-1}) \langle \Lambda_{-1} H, H \rangle +
\langle  (\Lambda_{-2}-\Lambda^{2}_{-1}) H, H \rangle +
\frac{1}{2} \langle \Lambda_{-1} H, H \rangle^2 \right]d\mu\,,
\end{array}\end{equation} where \begin{equation} d\mu = \pi^{-{m+2\choose
2}}\left|\det(2H'H'{}^*-|x|^2I)\right|\,e^{- \langle H, H \rangle}
\,dH'\,dx\,.\label{dmu}\end{equation}

Recalling \eqref{Aexp} and \eqref{L0}--\eqref{L2}, we see that $b_1$ is of the
form
$$b_1 = \sum c_{j\bar j' q \bar q'}  R_{j\bar j' q \bar q'}\;,$$ where
$c_{j\bar j' q \bar q'}$ is universal. Since $b_1$ is also invariant under the
unitary group, we must have
\begin{equation}\label{kappa1} b_1 = \be_1 \rho\;,\end{equation} where $\be_1$
is a universal constant (depending only on the dimension $m$ of $M$).  Similarly,
 $b_2$ is of the
form
$$b_2 = Q(R,R) +\ga_0\Delta\rho\;,$$ where $Q(R,R)$ is a universal quadratic
form in the curvature tensor $R$.  But $b_2$ is also $\U(m)$-invariant and
hence is a curvature invariant (of order 2).  Thus,
 \begin{equation} \label{kappa-two}b_2 =\ga_0
\Delta\rho+\ga_1\rho^2 +\ga_2|R|^2+\ga_3
|Ric|^2\;,\end{equation}
where the $\ga_k$ are universal constants depending only on $m$, which
completes the proof of Theorem
\ref{UNIVDIST}.

The proof of Theorem
\ref{univmorse} is exactly as above, except we integrate
over ${\bf S}_{m,q-m}$ instead of $\sym(\C,m)\times\C$ in the
computation of the expansion of $\kcal^\crit_{N,q}$.\qed

\subsection{Number of critical points: proof of Corollary \ref{EN}}

The coefficient $b_{1q}$ is of the form
\begin{equation}\label{b1q}b_{1q} = \beta_{1q}\,\rho\;,\end{equation}
by the above argument (or by the fact that $\rho$ is the only curvature
invariant of order 1).  Furthermore, it is well known (see, e.g.,
\cite[pp.~112--113]{Ko}) that for any \kahler metric $\om$ on $M$,
we have
\begin{equation}\label{chernforms}(\rho^2-|Ric|^2)\Om  = c_1(M,\om)^2\wedge
\om^{m-2}\,,\quad
   (|Ric|^2 -  |R|^2)\Om = [c_1(M,\om)^2-2c_2(M,\om)] \wedge
\om^{m-2}\;,\end{equation} where
$\Om=\frac1{4\pi^2m(m-1)}\,\om^m$. Therefore
 \begin{equation}\label{kappa2} b_{2q} =\ga_{0q}
\Delta\rho+(\ga_{1q}+\ga_{2q}+\ga_{3q})\rho^2 +\mbox{const.
}\frac{c_1(h)^2\wedge\om_h^{m-2}}{\om_h^m}+\mbox{const.
}\frac{c_2(h)\wedge\om_h^{m-2}}{\om_h^m}\;,\end{equation} where we
now write $c_j(h)=c_j(M,\om_h)$ for the $j$-th Chern form of the
\kahler metric $\om_h=-\frac i2\ddbar\log h$.

Integrating \eqref{kappa1} and \eqref{kappa2}, noting that $\int
\Delta \rho\,dV_h= \frac 2{m!}\int \ddbar \rho\wedge \om_h^{m-1} = 0$,
we obtain the asymptotic expansion of Corollary \ref{EN} with
$\be_{2q}=\ga_{1q}+\ga_{2q}+\ga_{3q}$.\qed

\subsection{Asymptotic expansions on Riemann surfaces: proof of Theorem \ref{critRS}}

On a compact Riemann surface $C$ of genus $g$, Corollary \ref{EN} says that
$\ncal^\crit_{N,q,h}$ has a universal expansion of the form
$$\ncal^\crit_{N,q,h} \sim \pi b_{0q}\,c_1(L)\,N + \pi\be_{1q}\,(2-2g)
+ \be_{2q}\int_C\rho^2\,\om_h\,N^{-1} +\cdots\,,$$ for $q=1,2$. There are
several ways to compute the constants.  A quick way to find
$b_0,\be_1,\be_2$, is to consider the case of $\CP^1$ with the Fubini-Study metric on
$L=\ocal(1)$. By an elementary computation in \cite{DSZ} (or by
\S\ref{s-exact}), we showed that for this case
\begin{eqnarray*}\ncal^\crit_{N,1,h}&=&
\frac{4(N-1)^2}{3N-2}\  =\ \frac 43\, N - \frac{16}{9} + \frac
4{27}\, N^{-1}\cdots \quad (\mbox{expected number of saddle points),}
\\ \ncal^\crit_{N,2,h}&=&\frac{N^2}{3N-2}\  =\ \frac 13\, N +
\frac{2}{9} + \frac 4{27}\, N^{-1} \cdots \quad (\mbox{expected number of
local maxima).}
\end{eqnarray*}
Therefore, $$\pi b_{01}=\frac 43,\quad \pi b_{02}=\frac 13,\quad
\pi\beta_{11}= -\frac 89,\quad \pi \beta_{02}=\frac 19\;.$$

To find $\beta_{21},\beta_{22}$, we note that $\int_{\CP^1}\om_\FS = \pi
c_1(L)=\pi$, where
$\om_\FS$ is the Fubini-Study \kahler form on $\CP^1$.  Furthermore, since
$c_1(\CP^1)=\frac 1\pi \int \rho \om_\FS=2$, we have $\rho\equiv 2$. (This
can be checked directly as follows: the \kahler potential
$K=\log(1+|z|^2) = |z|^2-\half |z|^4 +\cdots$, where $z$ is the affine
coordinate, and hence by \eqref{potexp},
$\rho(0)=R_{1\bar 1 1\bar1}(0)=2$.) Hence, $\int\rho^2\,\om_\FS=4\pi$, and
therefore
$$\be_{21}=\be_{22}= \frac 1{27\pi}\;,$$  which completes the proof of
Theorem \ref{critRS}.\qed

\section{Proof of Theorem \ref{ASYMPMIN}: evaluating the
coefficient $\be_{2q}(m)$}\label{s-coeff}

We have already shown that
\begin{itemize} \item \ $\int_M b_{1q} dV_h$ is topological;
\item \ $\int_M b_{2q} dV_h$ is the sum of a topological term plus a positive
multiple of $\int_M \rho^2_h dV_h$.
\end{itemize}

To complete the proof of Theorem \ref{ASYMPMIN} and  show that the
metric with asymptotically minimal $\ncal^\crit_{N,q,h}$ (and
$\ncal^\crit_{N,h}$) is the one for which $\om_h$ has minimal $\lcal^2$
norm of the scalar curvature, we must show that the $\be_{2q}$ are
positive.

The proof consists of finding an integral formula for $\be_{2q}$
and then transforming it to one that is amenable to (computer)
evaluation. We first summarize the key results:

\begin{lem}\label{K0}  In all dimensions,
\begin{equation} \label{I0} \be_{2q}(m)=  \frac
1{4\,\pi^{m+2\choose 2}}  \int_{{\bf S}_{m,q-m}'} \ga(H) \,
\left|\det(2HH^*-|x|^2I) \right| \,e^{ -\langle (H, x),(H, x)
\rangle}\, dH\, dx\;,\end{equation} where
$$\gamma(H) = \frac{1}{2} |H_{11}|^4 - 2 |H_{11}|^2 +1\,,$$ and
${\bf S}_{m,q-m}'$ is given by \eqref{Smk}.
\end{lem}

After a sequence of  manipulations as in the proof of Lemma
\ref{IZint}, the integral (\ref{I0}) will be rewritten in the
following form:
\begin{lem} \label{LI4}
$$\be_{2q}(m)= \frac {(-i)^{m(m-1)/2}} {4\,\pi^{2m}\prod_{j=1}^{m-1} j!}
\int_{Y_{2m-q}}\int_\R\cdots\int_\R\Delta(\la)\,\Delta(\xi)\,
\prod_{j=1}^m |\la_j|\;e^{i\langle
\la,\xi\rangle}\,\ical(\la,\xi)\, \,d\xi_1\cdots d\xi_m\,
d\la\;,$$ where \begin{equation}\label{Ilx} \ical(\la,\xi)=
\frac{F(D(\la))+\left[ \frac {4\sum_{j=1}^m \la_j} {m(m+1)(m+3)}
-\frac{2}{m+1}\right]\frac 1{\left(1 - \frac i2
\sum_j\xi_j\right)}+ \frac 2{(m+1)(m+3)\left(1 - \frac i2
\sum_j\xi_j\right)^2}} {\left(1 - \frac i2 \sum_j\xi_j\right)\prod_{j\le k}\left[ 1
+\frac i2 (\xi_j+\xi_k)\right] } \;.
\end{equation}
Here, $D(\lambda)$ is the diagonal matrix with diagonal entries
$\lambda = (\lambda_1, \dots, \lambda_m)$, $\Delta(\lambda) =
\Pi_{i < j} (\lambda_i - \lambda_j)$ is the Vandermonde
determinant and $$ F(P)= 1 - \frac {4\, \tr\, P }{m(m+1)}+\frac
{4(\tr\,P)^2+8\,\tr(P^2)} {m(m+1)(m+2)(m+3)}\;,$$for (Hermitian)
$m\times m$ matrices $P$.  The iterated $d\xi_j$ integrals are defined in the
distribution sense.
\end{lem}

The final step is the evaluation of $\beta_{2q}(m)$. Having
simplified the integral as far as we could, we complete the
computation  for the cases $m\le 4$ using Maple and find that it
is positive for these cases.
The resulting values of the  constants $\be_{2q}(m)$, $m\le 4$,
are given in \S \ref{s-beta2}.

\subsection{Proof of Lemma \ref{K0}}

We use the case of $M=\CP^1 \times E^{m-1}$ where $E$ is an
elliptic curve, and $L$ is the product of degree 1 line bundles on
the factors (with the Fubini-Study metric on $\ocal_{\CP^1}(1)$
and the flat metric on $E$). (The manifold $M$ is a homogeneous
space with respect to $SU(2) \times T^{2m-2}$, so the critical
point density is invariant and hence constant.)

Since $c_1(h)^2=c_2(h)=0$, it follows from \eqref{kappa2} that the
coefficient $b_{2q}$ of the expansion $N^{-m}\kcal^\crit_{N} (z) =
b_{0q}+b_{1q}N^{-1} +b_{2q}N^{-2} +O(N^{-3})$ is given by $b_{2q}=
\be_{2q}\rho^2$, and hence \begin{equation}\label{k0} \be_{2q}=\frac
1{\rho^2}\,b_{2q}=\frac 14\, b_{2q}\;.\end{equation}

The \szego kernel for $(M,L)$ is the product of the \szego kernels
on $\CP^1$ and $E^{m-1}$. Since the universal cover of $E^{m-1}$
is $\C^{m-1}$, the \szego kernel on $E^{m-1}$ is given by the
Heisenberg \szego kernel on $\C^{m-1}$ (see \cite[\S 1.3.2]{BSZ1})
modulo an $O(N^{-\infty})$ term, and we have:
$$ \Pi_{CP^1 \times E^{m-1}}(z,w) = \frac {(N+1)N^{m-1}}{\pi^m}
(1+z_1 \bar w_1)^N e^{N (z_2 \bar w_2+\cdots +z_m\bar w_m)
}\,e_L(z)\otimes \overline{e_L(w)} +O(N^{-\infty})
 \;.$$ As in \S \ref{s-exact}, we consider the
normalized \szego kernel
\begin{equation} \wt \Pi_N(z,w):=(1+z_1\bar w_1)^N
e^{Nz'\bar w'}, \qquad z'=(z_2,\dots,z_m),\ w'=(w_2,\dots,w_m)\;.
\end{equation}

We  have:
\begin{eqnarray*} \frac{\d\wt\Pi_N}{\d{z_1}} &=& N(1+z_1\bar
w_1)^{N-1}e^{N z' \bar w'}  \bar w_1\;,\\
\frac{\d\wt\Pi_N}{\d{z_\al }} &=& N(1+z_1\bar w_1)^{N} e^{N z'\bar w'} \bar w_\al \;,\\
\frac{\d^2\wt\Pi_N}{\d{z_1}\d\bar w_{1}} &=& \{N(1+z_1 \bar
w_1)^{N-1} +
N(N-1)(1+z_1 \bar w_1)^{N-2}  z_{1}\bar w_1\} e^{N z'\bar w'}\;,\\
\frac{\d^2\wt\Pi_N}{\d{z_\al } \d \bar w_{\al' }} &=&
\{N\de_{\al\al'} + N^2 z_{\al'} \bar w_\al  \}(1+z_1 \bar
w_1)^{N}e^{N z'\bar w'}\;,
\\
\frac{\d^2\wt\Pi_N}{\d{z_1}\d\bar w_\al} &=&
N^2 (1+z_1 \bar w_1)^{N-1} e^{N z_\al  \bar w_\al }  z_\al\bar w_1\;,\\
\frac{\d^2\wt\Pi_N}{\d{z_\al }\d\bar w_{1}} &=& \{N^2(1+z_1 \bar
w_1)^{N-1} z_1 \bar w_\al  \} e^{N z'\bar w'}\;,\\&&\hspace{3in}
2\le \al,\al' \le m\;.\end{eqnarray*} It suffices to compute the
density at $0$.  From the above, we have:
\begin{equation}\label{d4} \frac{\d^4\wt\Pi_N}{\d{z_j}\d z_q\d\bar
w_{j'}\d\bar w_{q'}} (0,0)=\left\{ \begin{array}{ll}2 N(N-1), \; &
j=q=j'=q'=1 \\
2N^2, & j=q=j'=q'>1\\
N^2, & j=j'\neq q'=q\\
\end{array} \right.\;.
\end{equation}

Recalling \eqref{AN}--\eqref{CN}, we then have:
\begin{eqnarray}\wt A_N(0) &=&
\left( \frac{\partial^2 \wt\Pi_N}{\partial z_j \partial
\bar{w}_{j'}}(0,0) \right) \ = \ NI\label{AN1}
\\ \wt B_N(0)&=&\left[\left(\tau_{jq}\,
\frac{\partial^3 \wt\Pi_N}{\partial z_j
\partial \bar{w}_{q'} \partial \bar{w}_{j'}} (0,0)\right)
\quad \left( N\, \frac{\partial \wt\Pi_N}{\partial
z_j}(0,0)\right) \right] \ =\ 0\,,\label{BN1}
\\[8pt]
\wt C_N(0)&=& \left[
\begin{array}{cc}\displaystyle
\left( \tau_{jq}\tau_{j'q'}\,\frac{\partial^4 \wt\Pi_N} {\partial
z_q
\partial z_j
\partial
\bar{w}_{q'}\partial \bar{w}_{j'}} (0,0) \right) &\displaystyle
\left( \tau_{jq}N\,\frac{\partial^2 \wt\Pi_N} {\partial z_j
\partial z_q}(0,0)
\right)
\\[8pt]
\displaystyle\left(\tau_{j'q'}N\,
 \frac{\partial^2 \wt\Pi_N}{\partial \bar{w}_{q'}\partial
\bar{w}_{j'}} (0,0) \right) & N^2\,\wt \Pi_N(0,0)
\end{array}
\right] \,,\label{CN1}\\[8pt] && \qquad\qquad \qquad 1\le j\le m\,,
1\le j\le q\le m\,, 1\le j'\le q'\le m\,.\nonumber
\end{eqnarray}

It follows that
\begin{equation}\label{La} \wt\Lambda_N(0) = \wt
C_N(0)=D(2N(N-1), \overbrace{2N^2,\dots,2N^2}^{(m-1)(m+2)/2},N^2)\
,
\end{equation}
i.e.\ the diagonal matrix with diagonal entries $2N(N-1)$, $2N^2$
repeated  $(m-1)(m+2)/2$ times, $N^2$.

We want to compute
$$\kcal^\crit_{N,q} (0) =
\frac{\pi^{-{m+2\choose 2}}}{\det \wt A_N(0) \det\wt \La_N(0)}
\int_{{\bf S}_{m,q-m}} \left|\det(HH^*-|x|^2I)\right|e^{
-{\langle \wt\La_N(0)^{-1}(H, x),(H, x) \rangle}}\,dH\,dx\,.
$$ Making the change of variables $H'\mapsto \sqrt2 \, N H',\
x\mapsto Nx$, the integral is transformed to
\begin{equation}\label{newI}\kcal^\crit_{N,q} (0) =
\frac{\pi^{-{m+2\choose 2}}N^m}{ \det\wh\La} \int_{{\bf S}_{m,q-m}'}
\left|\det(2HH^*-|x|^2I)\right|e^{ -{\langle
\wh\La^{-1}(H, x),(H, x) \rangle}}\,dH\,dx\,,
\end{equation} where
$$\wh \La = I-\frac 1N E\;,\qquad E=D(1,0,\dots,0)\;.$$
Therefore \begin{eqnarray*} N^{-m} \kcal^\crit_{N,q} (0) &=&
\pi^{-{m+2\choose 2}}\left(1+ \frac 1N + \frac 1{N^2} +
\cdots\right)\int_{{\bf S}_{m,q-m}'}
\left|\det(2HH{}^*-|x|^2I)\right|\\
&& \times \ \exp \left(- \|H\|^2-|x|^2-\frac 1N |H_{11}|^2 -
\frac 1{N^2}|H_{11}|^2 - \cdots\right) \,dH\,dx\\
&=& \left(1+ \frac 1N + \frac 1{N^2} + \cdots\right)\int
\exp\left(- \frac 1N |H_{11}|^2 - \frac 1{N^2}|H_{11}|^2 -
\cdots\right) \,d\wt\mu \\
&=& \int\left[1+\frac 1N(1-|H_{11}|^2)+ \frac 1 {N^2}\left(
1-2|H_{11}|^2 +\half |H_{11}|^4\right)\right]d\wt\mu + O\left(
\frac1{N^3}\right)\,,
\end{eqnarray*} where $$d\wt\mu= \pi^{-{m+2\choose
2}}\left|\det(2HH^*-|x|^2I) \right|
\,e^{ -\langle (H, x),(H, x) \rangle}\, dH\,
dx\;.$$
Therefore $$b_{2q}=\int_{{\bf S}_{m,q-m}'} \left( 1-2|H_{11}|^2 +\half
|H_{11}|^4\right)d\wt\mu\;,$$ and the desired formula then follows from
\eqref{k0}.\qed

\subsection{$\U(m)$ symmetries of the integral}

As an intermediate step between Lemmas \ref{K0} and \ref{LI4}, we
prove:

\begin{lem} \label{ICAL2}
\begin{equation}\label{I2}\be_{2q}(m)= \frac
1{4\,\pi^{m+2\choose 2}}  \int_{{\bf S}_{m,q-m}'} F(HH^*)\,
\left|\det(2HH^*-|x|^2I) \right| \,e^{ -\langle (H, x),(H, x)
\rangle}\, dH\, dx\;,\end{equation} where
\begin{equation}\label{F} F(P)= 1 - \frac {4\, \tr\, P }{m(m+1)}+\frac
{4(\tr\,P)^2+8\,\tr(P^2)} {m(m+1)(m+2)(m+3)}\;,\end{equation} for
(Hermitian) $m\times m$ matrices $P$.
\end{lem}

\begin{proof}

Since the change of variables $H\mapsto gHg^t$ ($g\in\U(m)$) is
unitary on $\sym(m,\C)$ (with respect to the Hilbert-Schmidt inner
product), we can make this change of variables in \eqref{I0},
and then integrate over
$g\in\U(m)$ to obtain
\begin{equation} \label{I1}\be_{2q}(m)= \frac
1{4\,\pi^{m+2\choose 2}} \int_{{\bf S}_{m,q-m}'}
\left(\int_{\U(m)}\ga(gHg^t)\,dg\right)
\,
\left|\det(2HH^*-|x|^2I) \right|
\,e^{ -\langle (H, x),(H, x) \rangle}\, dH\,
dx\;.\end{equation}

We now evaluate the integral $\int_{\U(m)}\ga(gHg^t)\,dg;$

\begin{claim}\label{G24} For $H\in \sym(m,\C)$, \begin{eqnarray}
\label{G2} \int_{\U(m)} |(gHg^t)_{11}|^2\,dg &=& \frac 2{m(m+1)} \tr
(HH^*)\;,\\
\label{G4}\int_{\U(m)} |(gHg^t)_{11}|^4\,dg &=&\frac
{8\, (\tr\, HH^*)^2 + 16
\,\tr(HH^*HH^*)}{m(m+1)(m+2)(m+3)}\;.\end{eqnarray}
\end{claim}

To prove the claim, we write $v=(v_1,\dots,v_m)=(g_{11},\dots,
g_{1m})$ so that $(gHg^t)_{11}=vHv^t$, and we replace
$\int_{\U(m)}dg$ with $\int_{ S^{2m-1}}d\nu(v)$, where $d\nu$ is
Haar probability measure on $ S^{2m-1}$. Next we recall that if
$p$ is a homogeneous polynomial of degree $2k$ on $\R^{2m}$,
\begin{equation}\label{stogauss} \int_{ S^{2m-1}}p(v)\,d\nu(v) = \frac
{(m-1)!} {(m-1+k )!} \int_{\R^{2m}}p(v)\,d\ga(v)\;,\qquad d\ga(v)=
\frac {1}{\pi^m}\,e^{-\|v\|^2}\,dv\;.\end{equation}

We easily see using Wick's formula that
\begin{eqnarray*}\int_{\C^m}|vHv^t|^2 \,d\ga &=& \sum_{j,k,j',k'}
H_{jk}\bar H_{j'k'}\int_{\C^m} v_jv_k\bar v_{j'}\bar v_{k'}\,d\ga \\
&=&
\sum_{j} |H_{jj}|^2\int_{\C^m} |v_j|^4\,d\ga  +2\sum_{j\neq k}
|H_{jk}|^2\int_{\C^m} |v_j|^2|v_k|^2\,d\ga\\&=& 2\,\tr (HH^*)\;.
\end{eqnarray*}
Formula \eqref{G2} then follows from \eqref{stogauss} with $k=2$.

Although the above approach can also be used to  verify \eqref{G4}, we
find it easier to use invariant theory, since the integral in
\eqref{G4} is a $\U(m)$-invariant function of $H\in\sym(m,\C)$, under
the $\U(m)$ action $H\mapsto gHg^t$.  Indeed, it is a
$\U(m)$-invariant Hermitian inner product on $S^2(\sym(m,\C))\approx
S^2(S^2 (\C^m))$.

 The action of $\U(m)$ on symmetric complex matrices
defines a representation equivalent to $S^2 (\C^m)$ where $\C^m$ is
the defining representation of $U(m)$. It is well known from
Schur-Weyl duality that $S^2(\C^m)$ is irreducible.
We then consider the $\U(m)$  representation $$S^2(S^2 (\C^m)) = \C\{
H_1 \otimes H_2 + H_2 \otimes H_1,\;\; H_1, H_2 \in S^2 (\C^m)\}, $$
with the diagonal action. Henceforth we put $$H_1 \cdot H_2 :=
\frac{1}{2} [H_1 \otimes H_2 + H_2 \otimes H_1]. $$ We then regard
$F(H)$ as the value on $H \otimes H$ of the quadratic form
$$Q( H_1 \cdot H_2) =
\int_{\U(m)}| \langle g H_1 g^t \cdot g H_2 g^t \,e_1 \otimes e_1,
e_1 \otimes e_1 \rangle|^2 dg. $$ This defines the Hermitian inner
product
$$\langle \langle H_1 \cdot H_2, H_2 \cdot H_4 \rangle \rangle =
\int_{\U(m)} \langle g H_1 g^t \cdot g H_2 g^t\,e_1 \otimes e_1, e_1
\otimes e_1 \rangle \overline{\langle g H_3 g^t \cdot g H_4 g^t\,
e_1 \otimes e_1, e_1 \otimes e_1 \rangle}  dg.
$$

We next recall that $S^2(S^2 (\C^m))$ decomposes into a direct sum of
two $\U(m)$ irreducibles, one corresponding to the Young diagram
$Y_1$  with $1$ row of four boxes and one corresponding to the
diagram $Y_2$ with $2$ rows each with two boxes. See for instance
Proposition 1 of \cite{Howe}. The Young projectors are
respectively,
$$\left\{ \begin{array}{l} P_{Y_1} (H \otimes H)_{i_1 i_2 i_3 i_4}
= \sum_{\sigma \in S_4} H_{i_{\sigma(1)} i_{\sigma(2)}}
H_{i_{\sigma(3)} i_{\sigma(4)}}\\ \\
P_{Y_2} (H \otimes H)_{i_1 i_2 i_3 i_4} = \sum_{\sigma \in S_2
\times S_2} (-1)^\sigma H_{i_{\sigma(1)} i_{\sigma(2)}}
H_{i_{\sigma(3)} i_{\sigma(4)}}.
\end{array} \right. $$
For $Y_2$ the $S_2 \times S_2$ permutes $1 \iff 3, 2 \iff 4. $

Since an irreducible $\U(m)$ representation has (up to scalar
multiples) a unique invariant inner product, it follows that
$$\langle \langle , \rangle \rangle = c_1 \langle, \rangle_{Y_1} +
c_2 \langle, \rangle_{Y_2}, $$ where $\langle, \rangle_{Y_j}$ are
the invariant inner products $$\langle A, B \rangle_{Y_j} = \tr\,
\Pi_{Y_j} (A) B^*$$ for the irreducibles corresponding to the
Young diagrams $Y_j$ as above.

We now calculate these inner
products on $H \otimes H$. We have
$$\left\{ \begin{array}{l} ||H \otimes H||_{Y_1}^2 = \sum_{\sigma \in S_4}
\sum_{i_1, i_2, i_3, i_4 = 1}^m H_{i_{\sigma(1)} i_{\sigma(2)}}
H_{i_{\sigma(3)} i_{\sigma(4)}} \bar{H}_{i_1 i_2} \bar{H}_{i_3 i_4}\\ \\
||H \otimes H||_{Y_2}^2  = \sum_{\sigma \in S_2 \times S_2}
\sum_{i_1, i_2, i_3, i_4 = 1}^m (-1)^\sigma H_{i_{\sigma(1)}
i_{\sigma(2)}} H_{i_{\sigma(3)} i_{\sigma(4)}} \bar{H}_{i_1 i_2}
\bar{H}_{i_3 i_4}.
\end{array} \right. $$
It is easy to see that each of these expressions is a linear
combination of the two quadratic forms
$$H \otimes H \mapsto \tr \{[H \otimes H] \circ [H^* \otimes
H^*]\},\;\;\; H\otimes H \mapsto [Tr H \circ H^*]^2. $$

Hence $$\int_{\U(m)} |(gHg^t)_{11}|^4\,dg = c_1\, (\tr\, HH^*)^2 + c_2
\,\tr(HH^*HH^*)\;.$$  To determine the constants $c_1,c_2$, it
suffices to consider the case where $H$ is diagonal.  Let
$s_1,\dots,s_m$ denote the eigenvalues of $H$.  Then by Wick's
formula we obtain
\begin{eqnarray*}\int_{\C^m}|vHv^t|^4 \,d\ga &=& \sum_{j,k,j',k'}
s_js_k\bar s_{j'}\bar s_{k'}\int_{\C^m} v_j^2v_k^2\bar v_{j'}^2\bar
v_{k'}^2\,d\ga
\\ &=&
\sum_{j} |s_j|^4\int_{\C^m} |v_j|^8\,d\ga  +2\sum_{j\neq k}
|s_j|^2 |s_k|^2\int_{\C^m} |v_j|^4|v_k|^4\,d\ga\\&=&
4!\sum_{j} |s_j|^4 +8\sum_{j\neq k}
|s_j|^2 |s_k|^2\
\\&=&
8\,(\tr\, HH^*)^2+16\,\tr\, (HH^*)^2\;.
\end{eqnarray*}
Formula \eqref{G4} now follows from \eqref{stogauss} with $k=4$.

Having proved the claim, the formula stated in Lemma \ref{ICAL2}
now follows from \eqref{I1} and Lemma \ref{G24}.
\end{proof}

\subsection{Proof of Lemma \ref{LI4}}

We  proceed exactly as in the proof of Lemma \ref{IZint}. We rewrite
the integral \eqref{I0} as
\begin{equation}\label{beta-ical}\be_{2q}(m)=
\frac 1{4\,\pi^m\,(2\pi)^{m^2}}\; \lim_{\ep'\to
0}\;\lim_{\ep\to
0}\ical_{\ep,\ep'}\;,\end{equation} where
\begin{eqnarray}
\ical_{\ep,\ep'} & = & \frac 1{\pi^{d_m}} \int_{{\mathcal H}_m}
\int_{\hcal_m(m-q)}\int_{\sym(m,\C)\times\C} \textstyle F(P+\half|x|^2I)\,
\left|\det(2P)\right|
 e^{i \langle \Xi, P - HH^* +\half |x|^2 I\rangle}\nonumber\\ && \times\
\exp\left(
-\tr  HH^*-|x|^2\right)\,
\exp\left({- \epsilon \tr \, \Xi\Xi^*-\ep' \tr  PP^* } \right)\, dH
\,dx
\,dP \,d
\Xi\;,\label{ical1}\\\hcal_m(m-q)&=& \{P\in \hcal_m:\mbox{index}\,P=m-q\}\;.
\nonumber\end{eqnarray}
Recall that $d_m= \dim_\C(\sym(m,\C)\times\C) = \half(m^2+m+2)$. As in
\S\ref{SIMPLIFY}, we note that absolute convergence is guaranteed by the
Gaussian factors in each  variable $(H, x, P,\Xi)$. Evaluating $\int
e^{i
\langle \Xi, P - HH^* +\half |x|^2\rangle}e^{-\ep\tr \, \Xi\Xi^*}d \Xi$
first, we obtain a dual Gaussian, which  approximates the delta function
$\de_{HH^* -\half |x|^2}(P)$. As $\epsilon \to 0$, the $d P$ integral
then yields the integrand at $P = H H^*-\half|x|^2I$; then letting
$\ep'\to 0$ we obtain the original integral stated in Lemma
\ref{ICAL2}.

Continuing as in \S\ref{SIMPLIFY}, we  conjugate $P$ to a diagonal matrix
$D(\lambda)$ with
$\lambda = (\lambda_1, \dots, \lambda_m)$ by an element $h \in
\U(m)$ and we replace $dP$ with $\Delta(\la)^2\,d\la\,dh$. Recalling
\eqref{IDENTITY}, we obtain:
\begin{eqnarray*}
\ical_{\ep,\ep'}
& = &  \frac{2^m\,c_m'}{\pi^{d_m}} \int_{\U(m)} \int_{{\mathcal H}_m}
\int_{Y'_{2m-q}}\int_{\sym(m,\C)\times\C}\Delta(\la)^2\,
\prod_{j=1}^m |\la_j|\,\textstyle F\left(D(\la)+\half |x|^2I\right) \\
&& \times\ e^{i\langle
\Xi,\, hD(\la)h^* +\half |x|^2I- HH^* \rangle}
e^{ -\left[\tr  HH^*+|x|^2+\epsilon
\tr
\Xi\Xi^* +\epsilon'  \sum \la_j^2\right]}\, dH \,dx \,d\la \,d
\Xi\,dh\;.
\end{eqnarray*}
Again using \eqref{IDENTITY} applied this time to $\Xi$, we obtain:
\begin{eqnarray*}
\ical_{\ep,\ep'}& = &  \frac{ 2^m(c'_m)^2}{\pi^{d_m}} \int_{\U(m)} \int_{\U(m)}
\int_{\R^{m}}\int_{Y'_{2m-q}}
\int_{\sym(m,\C)\times\C}\Delta(\la)^2\,\Delta(\xi)^2\,
\prod_{j=1}^m |\la_j|\,\textstyle F\left(D(\la)+\half |x|^2I\right)
\\&&\times\
e^{i \langle gD(\xi)g^*,\, hD(\la)h^* +\half |x|^2I- HH^*\rangle}
e^{ -\tr  HH^*-|x|^2-  \sum(\ep\xi_j^2
+\ep'\la_j^2)}
\, dH \,dx \,d\la\,d
\xi\,dh\,dg\;,\end{eqnarray*}
 We
then  transfer the conjugation  by $g$ to the right side of
the $\langle, \rangle$ in the first exponent and make the change of
variables
$h\mapsto gh, H\mapsto gHg^t$ to eliminate $g$ from the integrand:
\begin{eqnarray*} \ical_{\ep,\ep'}& = &   \frac{2^m(c'_m)^2 }{\pi^{d_m}}
\int_{\U(m)}
\int_{\R^{m}}\int_{Y'_{2m-q}}
\int_{\sym(m,\C)\times\C}\Delta(\la)^2\,\Delta(\xi)^2\, \prod_{j=1}^m
|\la_j|\;\textstyle F\left(D(\la)+\half |x|^2I\right)
\\&&
\times\
e^{i \langle D(\xi),\, hD(\la)h^* +\half |x|^2I- HH^*\rangle}
e^{ -\tr  HH^*-|x|^2- \sum(\ep\xi_j^2
+\ep'\la_j^2)}
\, dH \,dx \,d\la\,d
\xi\,dh\;.
\end{eqnarray*}

Next we substitute   the
Itzykson-Zuber-Harish-Chandra integral formula
\eqref{IZ} into
the above and expand
$$\det [e^{i
\xi_j \lambda_k}]_{jk} = \sum_{\sigma \in S_m} (-1)^\sigma\;
e^{i\langle  \xi, \sigma(\lambda)\rangle }, $$  obtaining a sum of
$m!$ integrals.  However, by making the change of variables
$\la'=\sigma(\la)$ and noting that $\De(\la')=  (-1)^\sigma
\De(\la)$, we see as before that the integrals of these terms are equal,
and so we obtain
\begin{eqnarray}
\ical_{\ep,\ep'}& = &    (-i)^{m(m-1)/2}\frac{c_m''}{\pi^{d_m}}
\int_{\R^{m}}\int_{Y'_{2m-q}}
\int_{\sym(m,\C)\times\C}\Delta(\la)\,\Delta(\xi)\,
\prod_{j=1}^m |\la_j|\;e^{i\langle
\la,\xi\rangle}\nonumber\\
&& \textstyle\times\ F\left(D(\la)+\half |x|^2I\right)
\exp\left(i \left\langle D(\xi),\,\half |x|^2I-HH^*
\right\rangle -\tr  HH^*-|x|^2\right)\nonumber\\&& \textstyle\times\
\exp\left( -\ep \sum\xi_j^2
-\ep'\sum\la_j^2\right)
\, dH \,dx \,d\la\,d
\xi\;.\label{ical3'}
\end{eqnarray} where
$$c_m''= \frac {2^{m^2}\,\pi^{m(m-1)}}
{\prod_{j=1}^m j!}\;.$$

 The phase
\begin{eqnarray}\Phi (H,x;\xi) &:=& i \left\langle D(\xi),
 \half |x|^2 I  -  HH^* \right\rangle
-\tr  HH^* -|x|^2\nonumber
\\&=& -\left[ \|H\|^2_\HS +i\sum_{j,k=1}^m
\xi_j|H_{jk}|^2 + \left(1 -\frac i2
\sum_j\xi_j\right)|x|^2\right]\nonumber
\\&=& -\left[\sum_{j\le
k}\left( 1 +\frac i2 (\xi_j+\xi_k)\right)|\wh H_{jk}|^2  +
 \left(1 - \frac i2 \sum_j\xi_j\right)|x|^2\right]\;.
 \label{phase}\end{eqnarray}
Thus,
\begin{equation}\label{ical4} \ical_{\ep,\ep'} =  (-i)^{m(m-1)/2}
c_m''
\int_{Y'_{2m-q}} \int_{\R^{m}}\Delta(\la)\,\Delta(\xi)\,
\prod_{j=1}^m |\la_j|\;e^{i\langle
\la,\xi\rangle}\,\ical(\la,\xi)\, e^ { -\ep \sum\xi_j^2
-\ep'\sum\la_j^2}
\,d\xi\, d\la\;,
\end{equation} where,
\begin{eqnarray*} \ical(\la,\xi)&=& \frac 1{\pi^{d_m}} \int_\C
\int_{\sym(m,\C)}F\left(D(\la)+\half |x|^2I\right)
e^{\Phi (H,x;\xi)} \,dH\,dx\\ &=&
\frac{1}{\prod_{j\le
k}\left( 1 +\frac i2 (\xi_j+\xi_k)\right)}
\int_\C F\left(D(\la)+\half |x|^2I\right)\, e^{-\left(1 - \frac i2
\sum_j\xi_j\right)|x|^2}\,dx\;.\end{eqnarray*}
To evaluate the $dx$ integral, we first expand the amplitude:
\begin{eqnarray*} F\left(D(\la)+\half |x|^2I\right)
 &=& F(D(\la))+\left[ \frac {4\sum_{j=1}^m \la_j} {m(m+1)(m+3)}
-\frac{2}{m+1}\right]|x|^2\\&&+ \frac 1{(m+1)(m+3)} |x|^4
\;,\end{eqnarray*} and then integrate to obtain \eqref{Ilx}.

To evaluate $\lim_{\ep,\ep'\to 0+}\ical_{\ep,\ep'}$,  we
first observe as in \S \ref{s-exact} that the map
$$(\ep_1,\dots,\ep_m)\mapsto
\int_{\R^{m}}\Delta(\xi)\, \;e^{i\langle
\la,\xi\rangle}\,\ical(\la,\xi)\, e^ { -\sum\ep_j \xi_j^2}
\,d\xi$$ is a continuous map from  $[0,+\infty)^m$
to the tempered distributions. Hence by \eqref{beta-ical} and \eqref{ical4}, we
have:
\begin{eqnarray}\be_{2q}(m)&=&  \frac
{(-i)^{m(m-1)/2}} {4\,\pi^{2m}\,\prod_{j=1}^m
j!}\;\lim_{\ep'\to 0^+}\;\lim_{\ep_1,\dots,\ep_m\to 0^+}
m!\int_{Y_{2m-q}}d\la\nonumber\\&&\times\  \int_{\R^{m}}
\Delta(\la)\,\Delta(\xi)\,
\prod_{j=1}^m |\la_j|\;e^{i\langle
\la,\xi\rangle}\,\ical(\la,\xi)\, e^ { - \sum\ep_j\xi_j^2
-\ep'\sum\la_j^2}
\,d\xi\;.\label{I3}\end{eqnarray} Letting $\ep_1\to 0,\dots,\ep_m\to
0,\ep'\to 0$ sequentially, we obtain the formula of Lemma \ref{LI4}.\qed

\subsection{Values of $\be_{2q}(m)$}\label{s-beta2} We use the
integral formula of  Lemma \ref{LI4} to compute the constants
$\be_{2q}(m)$. The $\xi_j$ integrals can be evaluated  using residues
as in \S \ref{residues}; the resulting $\la$ integrand is a
polynomial function of the $\la_j$ and $e^{\la_j}$. The integrals were
evaluated in dimensions $\le 4$ using 
Maple~10.\footnote{The Maple program used for this computation ({\tt
coefficient.mw}) is included in the source file of the arXiv posting
at \ {\tt http://arxiv.org/e-print/math/0406089}.}

In dimension 1, we reproduce the result from \cite{DSZ}:
$$\be_{21}(1)=\frac 1{3^3\cdot\pi}\;,\quad \be_{22}(1)=\frac 1{3^3\cdot\pi}\;.$$

In dimension 2, we have:
$$\be_{22}(2)=\frac 1{2^3\cdot 5\cdot\pi^2}\;,\quad \be_{23}(2)=\frac
{2^4}{3^4\cdot 5\cdot\pi^2}
\;,\quad \be_{24}(2)=\frac{47}{2^3\cdot 3^4\cdot 5\cdot\pi^2}\;.$$

In dimension 3, we have:
$$\be_{23}(3)=\frac {2^2}{5^3\cdot\pi^3}\;,\quad \be_{24}(3)=\frac
{11\cdot 23}{2^5\cdot 5^3\cdot\pi^3}
\;,\quad \be_{25}(3)=\frac{2^9\cdot 7}{3^6\cdot 5^3\cdot\pi^3}
\;,\quad \be_{26}(3)=\frac{23563}{2^5\cdot 3^6\cdot 5^3\cdot\pi^3}\;.$$

In dimension 4, we have:
$$\begin{array}{c}\displaystyle \be_{24}(4)=\frac {2^2}{3^2\cdot 7\cdot\pi^4}\;,\quad \be_{25}(4)=\frac
{2\cdot 3\cdot 41}{5^6\cdot 7\cdot\pi^4} \;,\quad
\be_{26}(4)=\frac{1056667}{2^6\cdot 3^2\cdot
5^6\cdot\pi^4}\;,\\[8pt]
\displaystyle  \be_{27}(4)=\frac{2^{15}\cdot 937}{3^8\cdot
5^6\cdot 7\cdot\pi^4} \;,\quad
\be_{28}(4)=\frac{267828299}{2^6\cdot 3^8\cdot 5^6\cdot
7\cdot\pi^4} \;.\end{array}$$

This
completes the proof of Theorem
\ref{ASYMPMIN} for $m\le 4$.  B. Baugher \cite{BB} used Mathematica to compute the 
constants
$\be_{2q}(5)$, and found that they are positive as well.\qed

\medskip
\begin{rem} The values of the coefficient $\beta_2(m)$ for the expected total number
of critical points are:
$$\begin{array}{c}\displaystyle\be_2(1) = \frac 2 {3^3\cdot \pi}\;, \quad
\be_2(2)=\frac{32}{405\pi^2} =\frac {2^5}{3^4\cdot 5\cdot\pi^2}
\;,\\[8pt]\displaystyle \be_2(3)= \frac{104}{729\pi^3} = \frac{2^3\cdot
13}{3^6\cdot\pi^3}\;,\quad \be_2(4)=\frac {17152}{45927\pi^4}=
\frac{2^8\cdot 67}{3^8\cdot 7\cdot\pi^4}\;.\end{array}$$\end{rem}

\section*{Appendix 1. \ Explicit formulas for $\CP^m$}

We give in this
appendix the precise formulas for the expected numbers ${\mathcal
N}^\crit_{N,q}(\CP^m)$ of critical points of  Morse index
$q$ of random sections of $H^0(\ocal(N), \CP^m)$ for $m\le 4$, as well as
the leading coefficients in dimensions $\le 6$. These formulas were obtained
by evaluating the integral in Proposition
\ref{numcrit} using
Maple~~10.\footnote{The Maple program used for this computation
({\tt exactprojective.mw}) is included in the source file of the
arXiv posting at \ {\tt http://arxiv.org/e-print/math/0406089}.}

 For $m=1$, we reproduce
the result from \cite{DSZ}:
$${\mathcal  N}^\crit_{N,1}(\CP^1)=
\frac{4(N-1)^2}{3N-2}\;,\quad
{\mathcal
N}^\crit_{N,2}(\CP^1)=\frac{N^2}{3N-2}\;;\quad \mbox{hence }\ {\mathcal
N}^\crit_{N}(\CP^1)=\frac{5N^2-8N+4}{3N-2}\;.
$$

For $m=2$, we obtain:
$$\textstyle{\mathcal  N}^\crit_{N,2}(\CP^2)=
{\frac {3\, \left( N-1 \right) ^{3}}{ \left( 2\,N-1 \right)
}}\;,\quad {\mathcal  N}^\crit_{N,3}(\CP^2)=
{\frac {16\, \left( N-1 \right) ^{3}{N}^{2}}{\left( 3\,N-2
 \right) ^{3}}}\;,\quad {\mathcal  N}^\crit_{N,4}(\CP^2)=
{\frac {{N}^{5} \left( 5\,N-4 \right) }{ \left( 3\,N-2
 \right) ^{3} \left( 2\,N-1 \right) }}\;.$$  Hence, the expected total number of
critical points is:
$$\ncal^\crit_N(\CP^2)=
{\frac {59\,{N}^{5}-231\,{N}^{4}+375\,{N}^{3}-310\,{N}^{2}+132\,N-24}{
 \left( 3\,N-2 \right) ^{3}}}\;.$$
To check the computation, we note that
$${\mathcal  N}^\crit_{N,2}(\CP^2)-{\mathcal  N}^\crit_{N,3}(\CP^2)+{\mathcal
N}^\crit_{N,4}(\CP^2) = N^2-3N+3= c_2(T^{*1,0}_{\CP^2}\otimes
\ocal(N))\;.$$

In the case  $m=3$, we obtain:
$$\begin{array}{rl}{\mathcal  N}^\crit_{N,3}(\CP^3)=
{\frac {8\, \left( N-1 \right) ^{4}}{ \left( 5\,N-2 \right) }}
\;,& {\mathcal  N}^\crit_{N,4}(\CP^3)= {\frac { \left( N-1 \right)
^{4}{N}^{2} \left( 63\,{N}^{2} -50\,N +10
 \right) }{ \left( 2\,N-1 \right) ^{4} \left( 5\,N-2
 \right) }}
\;,\\[8pt] {\mathcal  N}^\crit_{N,5}(\CP^3)=
{\frac {256\, \left( N-1 \right) ^{4}{N}^{5}}{ \left( 5\,N-2
 \right)  \left( 3\,N-2 \right) ^{5}}}\;,&
 {\mathcal  N}^\crit_{N,6}(\CP^3)=
{\frac {{N}^{9} \left(451\,{N}^{4}-1248\,{N} ^{3} +1280\,{N}^{2}
-576\,N+96\right) }{ \left( 2\,N-1 \right) ^{4} \left( 3\,N-2
 \right) ^{5} \left( 5\,N-2 \right) }}
\;.\end{array}$$ The expected total number of critical points is:
$$\textstyle\ncal^\crit_N(\CP^3)= \frac {637\,{N}^{8}-3978\,
{N}^{7}+11022\,{N}^{6}-17608\,{N}^{5}+17736\,{N}^{4}-11552\,
{N}^{3}+4768\,{N}^{2}-1152\,N+128} { \left( 3\,N-2 \right)
^{5}}\;.$$ To check the computation:
$$\sum_{q=3}^6{\mathcal  N}^\crit_{N,q}(\CP^3)= N^3-4N^2+6N-4 =
c_3(T^{*1,0}_{\CP^3}\otimes \ocal(N))\;.$$

For $m=4$, we obtain:
$$\begin{array}{l}{\mathcal  N}^\crit_{N,4}(\CP^4)=
{\frac{ 5\, \left( N-1 \right) ^{5}}{ \left( 3\,N-1
 \right) }}
\;,\quad {\mathcal  N}^\crit_{N,5}(\CP^4)= {\frac {16\, \left( N-1
\right) ^{5}{N}^{2} \left( 183\,{N}^{2}-120\,N +20 \right) }{
\left( 5\,N-2 \right) ^{5}}}
\;,\\[8pt] {\mathcal  N}^\crit_{N,6}(\CP^4)=
{\frac { \left( N-1 \right) ^{5}{N}^{5} \left(
396227\,{N}^{7}-1078546\,{N}^{6} +1261212\,{N}^{5}
-821326\,{N}^{4} +321695\,{N}^{3}
-75780\,{N}^{2}+9940\,N-560\right) }{ \left( 5\,N-2 \right) ^{5}
 \left( 2\,N-1 \right) ^{7} \left( 3\,N-1 \right) }}
\;,\\[8pt]
 {\mathcal  N}^\crit_{N,7}(\CP^4)=
{\frac { 4096\,\left( N-1 \right) ^{5}{N}^{9} \left(
109\,{N}^{2}-102 \,N+24 \right) }{ \left( 5\,N-2 \right) ^{5}
\left( 3\,N-2
 \right) ^{7}}}
\;,\\[8pt] {\mathcal  N}^\crit_{N,8}(\CP^4)={\frac {\al}{ \left( 5\,N-2
\right) ^{5} \left( 3\,N-2
 \right) ^{7} \left( 2\,N-1 \right) ^{7} \left( 3\,N-1 \right) }}\;,\\[8pt]
\scriptstyle \qquad \al= {N}^{14}
(14251551\,{N}^{10}-86984891\,{N}^{9}+237134546\,{N} ^{8} -
380216704\,{N}^{7}+397067360\,{N}^{6}-282219280\,{N}^{5}+138269792\,{N}^{4}\\
\scriptstyle\qquad\qquad -46114432\,{N}^{3}+
10020608\,{N}^{2}-1281280\,N+73216)\;.
\end{array}$$
The expected total number of critical points is:
\begin{eqnarray*}\ncal^\crit_N(\CP^4)&=& \scriptstyle(6571\,{N}^{11}-56373\,
{N}^{10}+221376\,{N}^{9}-524190\,{N}^{8}
+831075\,{N}^{7}-926382\,{N}^{6}+741276\,{N}^{5}-426392\,{N}^{4}\\&&
\scriptstyle + 173200\,{N}^{3}-47520\,{N}^{2}+8000\,N-640)/ (
3\,N-2) ^{7 } \;.\end{eqnarray*} Again, to check the computation:
$$\sum_{q=4}^8{\mathcal  N}^\crit_{N,q}(\CP^4)= N^4-5N^3+10N^2-10N+5=
c_4(T^{*1,0}_{\CP^4}\otimes \ocal(N))\;.$$

Formulas for dimensions 5 and 6 were also obtained using Maple, but they
are too lengthy to be informative.  Instead, we give below a
table of approximate numerical values of the leading coefficients
$n_{m+r}(m)=\frac{\pi^m}{m!}\,b_{0,m+r}(m)$ of
${\mathcal  N}^\crit_{N,m+r}(\CP^m)$, for
$r=0,1,\dots,m$, for
$m\le 6$:

\bigskip\begin{center}
\begin{tabular}{|c||l|l|l|l|l|l|l||l|}\hline
\multicolumn{9}{|c|}{Leading coefficients of $\,{\mathcal  N}^\crit_{N,m+r}(\CP^m)$}\\[3pt]
\hline\hline
$r$ &\quad\ 0 &\quad\ 1&\quad\ 2&\quad\ 3&\quad\ 4&\quad\ 5&\quad\ 6&\ \
total\\ 
\hline  \hline
$m=1$ & 1.33333 &
0.33333 &\multicolumn{5}{|c||}{}& 1.66667  \\
 \hline $m=2$ &1.5 & 0.59259& 0.09259&\multicolumn{4}{|c||}{}&
2.18519\\
\hline $m=3$ &1.6 & 0.78750& 0.21070& 0.02320&\multicolumn{3}{|c||}{}&
2.62140\\
\hline $m=4$ & 1.66667 & 0.93696& 0.33019& 0.06533 & 0.00543
&\multicolumn{2}{|c||}{}& 3.00457
\\
\hline $m=5$ & 1.71429 & 1.05448 & 0.44235 & 0.11939 & 0.01844& 0.00121
&& 3.35015\\
\hline $m=6$ &  1.75 &1.14903 &0.54457&0.17979 & 0.03884& 0.00486 & 0.00026
& 3.66734\\
\hline\end{tabular}
\end{center}

\bigskip

\section*{Appendix 2. \ Baugher's conjecture}

 In this appendix, we describe Baugher's
conjectured identity and its implications for the positivity of
$\beta_{2q}(m)$.
After calculating the integral of each of the individual terms in
Lemma \ref{LI4} using Mathematica~5, Baugher made  the
following conjecture \cite{BB}:

\renewcommand{\thetheo}{A.1}
\begin{conj} \label{C1}
$$\be_{2q}(m)= \frac {(-i)^{m(m-1)/2}} {4\,\pi^{2m}\prod_{j=1}^{m-1} j!}
\int_{Y_{2m-q}}\int_\R\cdots\int_\R\Delta(\la)\,\Delta(\xi)\,
\prod_{j=1}^m |\la_j|\;e^{i\langle \la,\xi\rangle}\,\jcal(\xi)\,
\,d\xi_1\cdots d\xi_m\, d\la\;,$$ where \begin{equation*}
\jcal(\xi)= \frac{4} {(m+1)(m+2)(m+3)\left(1 - \frac i2
\sum_j\xi_j\right)^2\prod_{j\le k}\left[ 1 +\frac i2
(\xi_j+\xi_k)\right] } \;.
\end{equation*}
\end{conj}

The conjecture has been verified in dimensions $m \leq 5$
\cite{BB}.  Further, Baugher evaluated this integral  by a lengthy
calculation:

\renewcommand{\thetheo}{A.2}
\begin{prop} \label{BB} \cite{BB} The integral on the right hand side
of Conjecture~\ref{C1}  equals
\begin{equation*} \frac{4}
{(m+1)(m+2)(m+3)} \int_{{\bf S}_{m,q-m}'} |x|^2\:
\left|\det(2HH^*-|x|^2I) \right| \,e^{ -\langle (H, x),(H, x)
\rangle}\, dH\, dx. \end{equation*}
 \end{prop}

The expression in Proposition~\ref{BB} is obviously
positive, and Conjecture~\ref{C1} thus implies $\be_{2q}(m) > 0$.

\end{document}